\documentclass[final,3p,times,number,sort&compress]{elsarticle}
\usepackage{amssymb}
\usepackage{amsmath}
\usepackage{amsfonts}
\usepackage{mathtools}
\usepackage{shuffle}
\usepackage{enumitem}
\usepackage{empheq}
\usepackage{wasysym}
\usepackage{verbatim}
\usepackage{graphicx}

\usepackage[
bookmarks=true,         
bookmarksnumbered=true, 
colorlinks=true, pdfstartview=FitV, linkcolor=blue, citecolor=blue,
urlcolor=blue]{hyperref}

\usepackage{microtype}
\usepackage{geometry}                
\usepackage{graphicx}
\usepackage{amssymb}
\usepackage{epstopdf}
\usepackage{amssymb}
\usepackage{amsmath}
\usepackage{amsfonts}
\usepackage{graphicx}
\usepackage{subfigure}
\usepackage{booktabs}
\usepackage{color}
\usepackage{array}
\usepackage{mathrsfs}
\usepackage{multirow}
\usepackage{float}
\usepackage{fancyhdr}

\newtheorem{theorem}{Theorem}[section]

\newtheorem{Theo}[theorem]{Theorem}
\newtheorem{Assu}[theorem]{Assumption}
\newtheorem{Rem}[theorem]{Remark}
\newtheorem{Lem}[theorem]{Lemma}

\newtheorem{example}[theorem]{Example}
\newtheorem{Exe}[theorem]{Example}
\newtheorem{problem}[theorem]{Problem}

\newtheorem{remark}[theorem]{Remark} 

\newcommand{\EE}{\mathbb{E}}

\newcommand{\RR}{\mathbb{R} }

\newcommand{\ga}{\gamma}

\newcommand{\si}{\sigma}

\let\Section=\section
\def\section{\setcounter{equation}{0}\Section}

\def\RR{\mathbb{R} }

\def\EE{\mathbb{E}}

\def\si{{\sigma}}

\def\vare{{\varepsilon}}

\journal{Communications in Nonlinear Science and Numerical Simulation}

\begin{document}
\begin{frontmatter}
\title {Positivity preserving  logarithmic Euler-Maruyama type scheme for stochastic differential equations }

\author {Yulian Yi$^a$}
\author {Yaozhong Hu$^b$}
\author {Jingjun Zhao$^a$ \corref{label:author}}
\cortext[cor1]{Corresponding author.\\
\emph{Email addresses}: \texttt{yulianyihit@163.com} (Y. Yi), ~\texttt{yaozhong@ualberta.ca} (Y. Hu),  ~\texttt{hit\_zjj@hit.edu.cn} (J. Zhao$^{*}$)}

\address{$^a$School of Mathematics, Harbin Institute of Technology, Harbin 150001, China}
\address{$^b$Department of Mathematics and Statistical Sciences, University of Alberta, Edmonton,T6G2G1,Canada }

\begin{abstract}
In this paper, we propose a class of explicit positivity preserving
numerical methods  
for general stochastic differential equations which have  positive solutions. Namely, all the numerical solutions are positive. Under some reasonable conditions, we obtain  the convergence and   the convergence rate results for   these methods. The main difficulty is to obtain  the strong convergence and  the convergence rate  for stochastic differential equations whose  coefficients
are of exponential growth. Some numerical experiments are provided to illustrate the theoretical results for our schemes.
\end{abstract}


\end{frontmatter}

\section{Introduction}
In this paper  we consider the numerical solution of the following 
scalar stochastic differential equation (SDE)
\begin{equation}\label{eq1}
\mathrm{d}x(t)=b(x(t))\mathrm{d}t+\sigma(x(t))\mathrm{d}W(t),
\end{equation}
where $b$ and $\sigma$ are locally Lipschitz continuous functions  and 
the initial data $x(0)=x_0>0$. 
In many practical problems the quantities (such as stock price in finance,
the population affected in epidemic and so forth)  that we model take positive values. When we use such a stochastic differential equation to model positive  quantity, first we need  to make sure that the solution to the equation is positive. 
This property of the equation can be guaranteed by the famous Feller's (non-explosion) criteria (see next section for more discussion).  Now the problem is that if the stochastic differential 
equation \eqref{eq1} has a positive solution, how to construct its numerical approximations so that they are all  positive? 

There have been some  works along this line.  On one hand, we know many existing numerical algorithms cannot maintain the property of positiveness for the approximate solutions. In fact,  in \cite{Kahl2008}, authors proved that the approximate solutions from the  Euler-Maruyama  method cannot remain positive for any (one dimensional) SDE. On the other hand,  there are several existing works on preservation of positivity. 
For instance, in \cite{Neuenkirch2014},  by using the Lamperti transformation, authors transformed the original SDE whose diffusion coefficient is strictly positive in a certain domain into the following SDE
\begin{equation*}
\mathrm{d}x(t)=F(x(t))\mathrm{d}t+c_{0}\mathrm{d}W(t)  .
\end{equation*}
 The transformed SDE was discretized by a drift implicit Euler method and was transformed  back  so that  an approximate solution is inside  the domain of the original solution under appropriate assumptions. The  convergence rate  of their scheme was proved to be   $1.0$  and their scheme  can be applied to practical models, such as the Cox-Ingersoll-Ross (CIR) model and the  constant elasticity of variance (CEV) model, etc. However,  sometimes an additional  nonlinear  equation needs  to be solved  in order to obtain the approximate solution, which is time consuming to some extent. A  class of balanced implicit methods which can preserve positivity under an appropriate  choice of step size and weight functions was shown  in \cite{Kahl2006}, but the convergence results were derived under the linear growth condition which is restrictive for some practical models. 
Since its importance in finance  the  following CIR model 

\begin{equation}\label{ciro}
\mathrm{d}x(t)=\kappa(\lambda-x(t))\mathrm{d}t+\theta\sqrt{\mid x(t)\mid}\mathrm{d}W(t), x(0)=x_{0}\, ,
\end{equation}
has received particular  attention.  
In the case of unattainable boundary ($2\kappa\lambda\geq\theta^{2}$), based on the Lamperti transformation and the semi-implicit Euler method, the following numerical scheme was proposed  in \cite{Alfonsi2005} for this model:
\begin{equation*}
x_{k+1}=\bigg( \frac{\sqrt{x_{k}}+\frac{\theta}{2}\Delta W_{k}}{2(1+\frac{\kappa}{2}\Delta)}+\sqrt{\frac{(\sqrt{x_{k}}+\frac{\theta}{2}\Delta W_{k})^{2}}{4(1+\frac{\kappa}{2}\Delta)^{2}}+\frac{(4\kappa\lambda-\theta^{2})\Delta}{8(1+\frac{\kappa}{2}\Delta)}}\bigg)^{2}\,, \quad 
k=0, 1, 2, \cdots.
\end{equation*}
It is obvious that the approximate solutions are all positive. 
 The convergence order of this scheme on some finite interval
 was obtained to be  $\frac{1}{2}$ in \cite{Dereich2012},  and  was improved to $1.0$ in \cite{Alfonsi2013}.
 This idea  was applied to SDEs  (driven by Brownian motions, see \cite{Neuenkirch2014}) and SDEs with poisson jumps (see \cite{Yang2017}) which take positive values.
Additionally, there are various modifications of explicit Euler-Maruyama method. For example, for the CIR process, the following explicit method 
 \begin{equation*}
 \begin{aligned}
 &\tilde{x}_{k+1}=\tilde{x}_{k}+\kappa\Delta(\lambda-(\tilde{x}_{k})^{+})+\theta\sqrt{(\tilde{x}_{k})^{+}}\Delta W_{k},~\tilde{x}_{0}=x_{0},\\
 &x_{k+1}=(\tilde{x}_{k+1})^{+},
 \end{aligned}
 \end{equation*}
 was proposed   in \cite{Lord2010}, and  in \cite{Cozma2020}   this scheme was proved to have strong convergent rate    $0.5$ in $L^{p}$ sense when  $\frac{2\kappa\lambda}{\theta^{2}}>3$, $1\leq p<\lambda-1$.  A symmetrized Euler method 
  \begin{equation*}
 \begin{aligned}
 &x_{k+1}=\mid x_{k}+\kappa\Delta(\lambda-x_{k})+\theta\sqrt{x_{k}}\Delta W_{k}\mid,
 \end{aligned}
 \end{equation*}
 with $L^{p}$ convergence rate  of $0.5$ was  obtained  in \cite{Berkaoui2008} under a more  restricted condition on  the 
 parameters 
 \begin{equation*}
 \frac{2\kappa\lambda}{\theta^{2}}>1+\sqrt{8}\max\{ \frac{\sqrt{\kappa}}{\theta}\sqrt{16p-1}, 16p-2\}.
 \end{equation*}
  Some other works can be found  in \cite{Moro2007, Halidias2014,Higham2013} and references therein.

However, for general stochastic differential equations  not much study  to construct positive  numerical   approximations   has  been done. In this work we shall construct   positive numerical  approximations for general equations.  Namely, once the equation \eqref{eq1} satisfies the Feller's condition so that its solution is positive,  we shall  construct its  numerical approximations which are all positive themselves.  More precisely, we shall  propose an explicit logarithmic Euler-Maruyama type scheme such that it can keep positive unconditionally.


This paper is organized as follows. In Section 2, we give a general idea of our method. In Section 3, we consider  the almost sure convergence of a class of positivity preserving numerical methods obtained by combining the logarithmic transformation and some classical stochastic methods, such as the Euler-Maruyama method and the stochastic theta method. In Section 4, we prove the exponentially integrability of both the exact solution and numerical solution from the logarithmic truncated Euler-Maruyama method (logTE). Moreover,  we obtain  the strong convergence results of the logTE method under Assumption \ref{as1} with $m<2$. Since there are    models (for example, the CIR process) which do   not satisfy this assumption, we weaken Assumption \ref{as1} to Assumption \ref{as4.2}, and prove that the logTE numerical approximation is also strongly convergent. Next, in Section 5, we analyze the convergence rate of the logTE method under a slightly stronger Assumption \ref{as2r}. Some numerical experiments are carried out in Section 6  to illustrate our theoretical conclusions for some popular models. Finally, we make a brief  conclusion in Section 7.
  
\section{General idea}
Since we consider only positive solution to the equation \eqref{eq1}, we  assume that the coefficients $b$ and $\si$ are real valued functions defined on the positive real line $\RR_+
:=(0, \infty)$  and we also  assume they  are  locally Lipschitz continuous.   
From these coefficients,  we introduce  the scale function
\begin{equation*}
p(x)=\int\limits_{c}^{x}\exp\left\{ -2\int\limits_{c}^{\xi}\frac{b(u)}{\sigma^{2}(u)}\mathrm{d}u\right\}\mathrm{d}\xi \,, \ \ x\in \RR_+,
\end{equation*}
for some fixed $c\in  \RR_+$ and define  
\begin{equation*}
v(x)=\int\limits_{c}^{x}p'(y)\int\limits_{c}^{y}\frac{2\mathrm{d}z}{p'(z)\sigma^{2}(z)}\mathrm{d}y \,, \ \ x\in \RR_+\,. 
\end{equation*} 
\begin{Assu}\label{as0} We assume that the following Feller's condition holds true:
\begin{equation}\label{eq4}
\lim\limits_{x\rightarrow0{+}}v(x)=+\infty
\quad  {\rm and}\quad    \lim\limits_{x\rightarrow+\infty}v(x)=+\infty\,. 
\end{equation} 
\end{Assu}
According to Feller's criteria (see \cite{karatzas1991brownian}) 
the condition \eqref{eq4} is equivalent to the fact  that  the SDE \eqref{eq1}  has a unique strong   solution $x(t), t\geq0$,
and its sample path is contained in $\RR_+$, i.e.
\begin{equation}\label{eq3}
\mathbb{P}\big(x(t)\in(0,\infty),t\geq0\big)=1.
\end{equation}
Now assume the condition \eqref{eq4} holds true, the solution is then 
positive for all time $t\ge 0$. How we construct numerical approximations so that they are also positive?
 In the introduction we briefly mentioned some   earlier works  about positivity preserving numerical schemes. In this work, we want to find   general  numerical schemes  to
construct positive approximate solutions under the Feller's condition criteria
\eqref{eq4}.  Our general idea is first to take  the following logarithmic 
transformation
\begin{equation*}
y=\phi(x)=\ln x,\qquad x\in \RR_+\,. 
\end{equation*}
The application of the  It\^o  formula with $y(t)=\ln x(t)$  yields  the following transformed stochastic differential equation for $y(t)$: 
\begin{equation}\label{eq5}
\begin{aligned}
\mathrm{d}y(t)=f(y(t))\mathrm{d}t+g(y(t))\mathrm{d}W(t),\quad y(0)=y_{0}=\ln x_{0},
\end{aligned}
\end{equation}
where 
\begin{equation}
f(y)=e^{-y}b(e^{y})-\frac{1}{2}e^{-2y}\sigma^{2}(e^{y})\,,\quad g(y)=e^{-y}\sigma(e^{y})\,. 
\end{equation}
Then we find  a numerical solution $\tilde y(t)$ for \eqref{eq5}. 
And $\tilde x(t)=e^{\tilde y(t)}$ will be a positive approximation of the original solution $x(t)$.   To realize this idea  the main difficulty
is to approximate \eqref{eq5} since the coefficients are  usually no longer Lipschitzian
even when the original ones are. But first, let us state the following result about \eqref{eq5}.  
\begin{Theo}\label{2.1}
If equation \eqref{eq1} satisfies Assumption \ref{as0}, then  SDE \eqref{eq5} admits an unique nonexplosive solution  in $\mathbb{R}$, i.e.
\begin{equation*}
\mathbb{P}(y(t)\in(-\infty,\infty),t\geq0)=1.
\end{equation*}
\end{Theo}
\textbf{Proof.}
Denote 
\begin{equation*}
\tilde{p}(y)=\int\limits_{c}^{y}\exp\left\{  -2\int\limits_{c}^{\xi}\frac{f(u)}{g^{2}(u)}\mathrm{d}u\right\}\mathrm{d}\xi  ,
\end{equation*}
and 
\begin{equation*}
\tilde{v}(y)=\int\limits_{c}^{y}\tilde{p}'(u)\int\limits_{c}^{u}\frac{2\mathrm{d}z}{\tilde{p}'(z)g^{2}(z)}\mathrm{d}u\,. 
\end{equation*}
Then the  nonexplosiveness of the solution $y(t)$  can still be explained by the Feller's condition, i.e.
\begin{equation*}
\lim\limits_{x\rightarrow-\infty}\tilde{v}(x)=\infty 
\quad  {\rm and}\quad \lim\limits_{x\rightarrow\infty}\tilde{v}(x)=\infty.
\end{equation*}
Since $y=\ln x$, integration by substitution yields
 \begin{equation*}
\lim\limits_{x\rightarrow-\infty}\tilde{v}(x)=\infty \Leftrightarrow \lim\limits_{x\rightarrow0_{+}}v(x)=\infty ,
\end{equation*}
and 
 \begin{equation*}  \lim\limits_{x\rightarrow\infty}\tilde{v}(x)=\infty  \Leftrightarrow\lim\limits_{x\rightarrow\infty}v(x)=\infty\,. 
\end{equation*}
Hence, we can derive that  the solution of SDE \eqref{eq5} will stay in $(-\infty,\infty)$ with probability 1. Namely, with probability 1, the solution $y(t)$ to \eqref{eq5} has no finite time explosion.\qed

Since the SDE \eqref{eq5} takes value in the whole space, we do not need to worry  whether the numerical approximation is contained in $\RR_+$ (or a sub-domain of $\RR$) or not. However,  after   the logarithmic transformation, the coefficients of the transformed equation 
\eqref{eq5} are usually of super-linear growth or even of exponential growth.  This brings new difficulty  to the construction of numerical scheme. Take the CIR model   \eqref{ciro} as an illustrating  example,  which was commonly used in financial and biological mathematics. Using the Feller condition, we can derive that when $\kappa\lambda\geq\frac{\theta^{2}}{2}$, $ x_{0}>0$, the CIR process has a (unique)  positive solution, see \cite{Andersen2007,Dereich2012}. The use of the logarithmic transformation  $y(t)=\ln x(t)$ will give 
\begin{equation}\label{cirt}
\mathrm{d}y(t)=[(\kappa\lambda-\frac{\theta^{2}}{2})e^{-y(t)}-\kappa]\mathrm{d}t+\theta e^{-\frac{1}{2}y(t)}\mathrm{d}W(t),~~y(0)=\ln x_{0}.
\end{equation}
The coefficients of the above equation increase exponentially
(as $y$ goes to $-\infty$), so the classical strong convergence theories based on the global Lipschitz condition or polynomial growth condition are no longer applicable to the above equation.

\section{Almost sure convergence }
First, let us introduce some known conclusion about the almost sure convergence of some classical stochastic methods. To simplify the presentation, we fixed a finite time interval $[0, T]$ for any given fixed $T>0$.  Consider the partition $\pi: 0=t_0<t_1<\cdots<t_N=T$,
where $N$ is a positive integer and denote $\Delta=T/N$ as the uniform step size of the interval $[0,T]$.

There are several papers considered the almost sure convergence of stochastic methods for SDE with non-global Lipschitzian coefficients, see \cite{Gyongy1998, Jentzen2009} and references therein. In \cite{Fleury2005}, author  proved the almost sure convergence of a family of numerical  methods, which   cover a number of   classical stochastic schemes. The methods has the following form
\begin{equation}\label{aeq5}
\left\{
\begin{aligned}
&y_{\Delta}^{\pi}(t_{i+1})=\Delta f_{\Delta}^{\pi}(y_{\Delta}^{\pi}(t_{i}))+g_{\Delta}^{\pi}(y_{\Delta}^{\pi}(t_{i}))\delta W_{i}+\varphi_{\Delta}^{\pi}(y_{\Delta}^{\pi}(t_{i}))\delta Z_{\Delta,i},\quad i=0,1,\cdots,N-1,\\
&y_{\Delta}^{\pi}(0)=y_{\Delta,0}^{\pi}\,,
\end{aligned}
\right.
\end{equation}
where  $\delta U_{i}=U(t_{i+1})-U(t_{i})$, $U=W \text {or } Z$, and other notations are explained  in the following assumption. 
\begin{Assu}\label{3aa}
(1) The coefficients $f_{\Delta}^{\pi}$, $g_{\Delta}^{\pi}$ and $\varphi_{\Delta}^{\pi}$ are locally bounded for $\Delta$ small enough:  $\forall R>0$, there exists a  $\Delta^{*}(R)$ such that for any $0<\Delta\leq \Delta^{*}(R)$
\begin{equation*}
\mid f_{\Delta}^{\pi}(x)\mid\vee\mid g_{\Delta}^{\pi}(x)\mid\vee\mid \varphi_{\Delta}^{\pi}(x)\mid\leq C_{R},\quad \forall 
\    \mid x\mid\leq R\, ,
\end{equation*}
for some $C_R\in (0, \infty)$. 

\noindent (2) $f_{\Delta}^{\pi}$, $g_{\Delta}^{\pi}$   converge quickly enough to $f$ and $g$   uniformly on any compact sets: for any  $R>0$, there exists a  $\Delta^{o}(R)$ such that, for any $0<\Delta<\Delta^{o}(R)$, 
\begin{equation*}
\mid f_{\Delta}^{\pi}(x)-f(x)\mid\vee\mid g_{\Delta}^{\pi}(x)-g(x)\mid\leq C_{R}\sqrt{\Delta},  \quad \forall 
\    \mid x\mid\leq R\,. 
\end{equation*} 

\noindent (3)  For any arbitrarily small $\Delta>0\in (0, 1]$,  $(Z_{\Delta}(t), t\in [0, T])$ is a  continuous $L^{2}$ martingale starting  at $0$ and   its quadratic variation process  $A_{\Delta}=[Z_{\Delta}]$ satisfies 
\begin{equation*}
\mathrm{d}A_{\Delta}(t)=h_{\Delta}(t)\mathrm{d}t, 
\end{equation*}
with 
\begin{equation*}
\sup_{u\in [0, T]}\EE \mid h_{\Delta}(u)\mid^{p}\leq C_{p, T}\Delta^{\frac{p}{2}},\quad \forall  \Delta\in(0,1] ,
\end{equation*}
for  any  positive integer  $p$.
 
\noindent  (4) For any positive integer $p$,  it holds that
\[
\EE\mid y_{\Delta,0}^{\pi}-y_{0}\mid^{p}\leq C_p \Delta^{\frac{p}{2}}
,\quad \forall  \Delta\in(0,1] .
\]
\end{Assu}

It is proved in   \cite[Theorem 2]{Fleury2005}  that scheme \eqref{3aa} converges almost surely.  We restate it as the following lemma.
\begin{Lem}\label{ascf} (\cite{Fleury2005})
If   $f, g$ satisfy the local Lipschitz condition and solution of \eqref{eq5} does not explode and if  Assumption \ref{3aa} holds,  then for any $\gamma\in(0,\frac{1}{2})$, there exists a random variable $\eta=\eta_{\gamma,T}$ such that
\begin{equation*}
\sup\limits_{i\in\{0,1,\cdots,N\}}\mid y_{\Delta}^{\pi}(t_{i})-y(t_{i})\mid\leq \eta \Delta^{\gamma},  \quad  \text a.s.
\end{equation*}
\end{Lem}
Then, we can get the following theorem.
\begin{Theo}
If Assumption \ref{as0} holds and the  coefficients of approximation process \eqref{aeq5} satisfy Assumption \ref{3aa}, then $x_{\Delta}^\pi (t_{i})=e^{y_{\Delta}^{\pi}(t_{i})}$ converges almost surely to
the solution of \eqref{eq1}.
More precisely,   for any $\gamma\in(0,\frac{1}{2})$, there exists a finite random variable $\tilde{\eta}_{\ga, T}$ such that
\begin{equation}
\sup\limits_{i\in\{o,1,\cdots,N\}}\mid x_{\Delta}^\pi(t_{i})-x(t_{i})\mid\leq \tilde{\eta}_{\ga, T} \Delta^{\gamma}, \quad \text a.s.\,. \label{e.3.4}
\end{equation}
\end{Theo}
\textbf{Proof.}
The conclusion of almost sure convergence for this theorem  can be obtained by combining the mean value theorem and Lemma \ref{ascf}. 
\qed

\begin{Exe}
Consider the Euler-Maruyama scheme
\begin{equation}\label{EM3}
y_{\Delta}^{E}(t_{i+1})=y_{\Delta}^{E}(t_{i})+\Delta f(y_{\Delta}^{E}(t_{i}))+g(y_{\Delta}^{E}(t_{i}))\delta W_{i},~~ y_{\Delta}^{E}(0)=y_{0}.
\end{equation}
 Let $f_{\Delta}^{\pi}=f$, $g_{\Delta}^{\pi}=g$, $\varphi_{\Delta}^{\pi}=0$, $Z_{\Delta}=0$, $y_{\Delta,0}^{\pi}=y_{0}$. Then obviously Assumption \ref{3aa} is satisfied. 
 \end{Exe}

\begin{Exe}
Consider the stochastic theta method
\begin{equation*}
y_{\Delta}^{IE}(t_{i+1})=y_{\Delta}^{IE}(t_{i})+\theta\Delta f(y_{\Delta}^{IE}(t_{i}))+(1-\theta)\Delta f(y_{\Delta}^{IE}(t_{i+1}))+g(y_{\Delta}^{IE}(t_{i}))\delta W_{i}.
\end{equation*}
Denote
\begin{equation*}
\Phi_{\Delta}(x)=x-(1-\theta)\Delta f(x).
\end{equation*}
According to analysis of \cite{Fleury2005}, we know that $\Phi_{\Delta}(x)$ is locally Lipschitz if $f$ is.  As a result, for any $R>0$, $\Phi_{\Delta}$ can be inverted on a ball $B_{R}=\{x: \mid x\mid \leq R\}$ for sufficiently small $\Delta$, and $\Phi^{-1}_{\Delta}$ is Lipschitz continuous on $\Phi_{\Delta}(B_{R})$  and 
\begin{equation*}
\Phi_{\Delta}(y_{\Delta}^{IE}(t_{i+1}))=y_{\Delta}^{IE}(t_{i})+\theta \Delta f(y_{\Delta}^{IE}(t_{i}))+g(y_{\Delta}^{IE}(t_{i}))\delta W_{i}.
\end{equation*}
Let $Y_{\Delta}^{\pi}=\Phi_{\Delta}(y_{\Delta}^{IE})$. If $\Delta$ is sufficiently small, we  have 
\begin{equation*}
Y_{\Delta}^{\pi}(t_{i+1})=Y_{\Delta}^{\pi}(t_{i})+\Delta f(\Phi^{-1}_{\Delta}(Y_{\Delta}^{\pi}(t_{i})))+g(\Phi^{-1}_{\Delta}(Y_{\Delta}^{\pi}(t_{i})))\delta W_{i}.
\end{equation*} 
Take $f_{\Delta}^{\pi}=f\circ\Phi_{\Delta}^{-1}$, $g_{\Delta}^{\pi}=g\circ\Phi_{\Delta}^{-1}$, $\varphi_{\Delta}^{\pi}=0$, $Z_{\Delta}=0$, $y_{\Delta,0}^{\pi}=\Phi_{\Delta}(y_{0})$. Then Assumption \ref{3aa} is satisfied, so $Y_{\Delta}^{\pi}(t_{i})\rightarrow y(t_{i})$, a.s. with order $\gamma$  for any $\ga\in (0, 1/2)$. Besides, we also  have $y_{\Delta}^{IE}(t_{i})\rightarrow Y_{\Delta}^{\pi}(t_{i})$, as $\Delta$ goes to $0$. Therefore, the stochastic theta approximation $y_{\Delta}^{IE}(t_{i})$ will converge almost surely to $y(t_{i})$, and $x_{\Delta}(t_{i})=e^{y_{\Delta}^{IE}(t_{i})}$ will converge almost surely to exact solution $x(t_{i})$ of the original  SDE \eqref{eq1}. 
\end{Exe} 

\begin{remark} The random constant $\tilde{\eta}_{\gamma, T}$ in inequality \eqref{e.3.4} may not be integrable. The path-wise convergence cannot imply the strong convergence or weak convergence. Therefore, in the next section, we shall study another kind of convergence: the convergence in $L^{p}$ sense.
\end{remark}

\section{ Strong convergence analysis}
In the precedent  section, we analyze the almost sure convergence of some numerical schemes. Path-wise convergence plays an important role in analysis of random dynamical systems \cite{Arnold1998}, while sometimes we need more. For example, in mathematical finance, the pricing of complex derivative products  is one of the main tasks, which means we have to compute $\EE \tilde{\varphi}(x)$, where $\tilde{\varphi}: C([0,T];\mathbb{R})\rightarrow\mathbb{R}$ represents the discounted payoff of the derivative. Weak convergence is usually enough for pricing purposes, but to minimize the computational complexity strong convergence results are required when using multilevel Monte Carlo schemes \cite{Giles2008}. But the integrability of random constants in inequality  \eqref{e.3.4} is still unknown, and it seems hard to analyze these constants. If $\tilde{\varphi}$ is continuous and bounded, then $\EE\tilde{\varphi}(x_{\Delta}(t))\rightarrow\EE \tilde{\varphi}(x(t))$ as $\Delta\rightarrow 0$. However,  since the p-th power function needed for the $L^{p}$ convergence is not bounded, the path-wise convergence results cannot imply strong or weak convergence. For example,   in \cite{Hutzenthaler2011Strong}, in terms of the Euler scheme \eqref{EM3}, authors prove the following lemma.
 \begin{Lem}(\cite{Hutzenthaler2011Strong})
Let $f, g: \mathbb{R}\rightarrow\mathbb{R}$, $g(y_{0})\neq 0$. There exist constants $C\geq1, \tilde{\alpha}>\tilde{\beta}>1$  such that
\begin{equation*}
\max\big( \mid f(y)\mid, \mid g(y)\mid \big)\geq \frac{1}{C}\mid y\mid ^{\tilde{\beta}},~ ~\min \big(\mid f(y)\mid, \mid g(y)\mid \big)\leq C\mid y\mid ^{\tilde{\alpha}},
\end{equation*}
for all  $\mid y\mid\geq C$. 
Then, there exists a constant $c\in(1,\infty)$ and a sequence of nonempty events $\Omega_{N}\in\mathcal{F}$, $N\in\mathbb{N}$, with $P[\Omega_{N}]\geq c^{(-N^{c})} $ and $\mid y_{\Delta}^{E}(\omega,T)\mid\geq 2^{(\tilde{\alpha}^{(N-1)})}$ for all $\omega\in\Omega_{N}$ and all $N\in\mathbb{N}$. Moreover, if exact solution of  the SDE \eqref{eq5} satisfies
\begin{equation*}
\sup_{ t\in[0,T]} \EE\mid y(t)\mid^{p}<\infty, 
\end{equation*}
for some $p\in[1,\infty)$, then
\begin{equation*}
\lim_{N\rightarrow\infty} \EE\mid y(T)-y_{\Delta}^{E}(T)\mid^{p}=\infty, ~~\lim_{N\rightarrow\infty} \EE \mid y(T)\mid ^{p}-\EE\mid y_{\Delta}^{E}(T)\mid^{p} =\infty. 
\end{equation*}
\end{Lem}

We already know that the Euler-Maruyama approximation  converges  almost surely to exact solution if coefficients are Local Lipschitzian. The above lemma states that the Euler-Maruyama approximation may  diverge in both strong and weak $L^p$ sense.

Since we are concerned with the $L^p$ convergence, we need to deal with the exponential integrability of the numerical solutions from our transformation.  Therefore, the existing  strong convergence conclusions (see, \cite{zong2015,xiao2016} and references there in) are not applicable to \eqref{cirt} because of its exponentially growing coefficients and because of the additional requirement of the exponential integrability.  We shall focus on the Euler-Maruyama method, which is favored in practice because of its simplicity. Even for this simplest scheme the exponential growth causes significant difficulty.

Hence,  we need to  introduce some additional structure so that  the Euler-Maruyama method   will converge to solutions of SDEs with exponentially growing coefficients. 
Our idea is  inspired by \cite{Mao2015,Mao2016}.
First, let us make some assumptions on the coefficients of SDE \eqref{eq5}.
\begin{Assu}\label{as1}
We assume that the coefficients of SDE \eqref{eq5} are locally Lipschitzian and satisfy the following condition:
for any   $m\geq1$, there exists a constant $K
=K_{  m}>0$  (independent of  $p$ and $y$)    such that 
\begin{equation}\label{moas}
yf(y)+\frac{p-1}{2}g^{2}(y)\leq K p^{m}\,, \quad \forall \ p>2 ,  \  y\in\mathbb{R}\,. 
\end{equation}
\end{Assu}


Let $l>0$ be an arbitrarily large positive number.  Consider the following stochastic differential equation: 
\begin{equation}\label{eq5'}
\mathrm{d}y_{l}(t)=f_{l}(y_{l}(t))\mathrm{d}t+g_{l}(y_{l}(t))\mathrm{d}W(t),
\end{equation}
where 
\begin{equation*}
f_{l}(x)=f\big((\mid x\mid\wedge l)\frac{x}{\mid x\mid}\big)\quad {\rm and}\quad   g_{l}(x)=g\big((\mid x\mid\wedge l)\frac{x}{\mid x\mid}\big).
\end{equation*}
Since coefficients of SDE \eqref{eq5} satisfy condition \eqref{moas}, and coefficients of SDE \eqref{eq5'} satisfy the global Lipschitz condition, so there exists a unique solution $y_{l}(t)$, which  converges to the exact solution y(t) of SDE \eqref{eq5} in probability
as $l$ goes to infinity.  If we apply the Euler-Maruyama method to \eqref{eq5'},  we obtain a numerical approximation $y_{\Delta, l}(t)$ to \eqref{eq5'} which will strongly converge to $y_{l}(t)$ in $L^{p}$ ($\forall p>0$) sense. Therefore,  for each $\Delta$, there is a  $l_{\Delta}$ such that $y_{\Delta, l_{\Delta}}(t)$ will converge to $y(t)$ in probability. 
  We will prove that under Assumption \ref{moas}, 
   both the exact solution and the numerical solution are exponentially integrable and  like in  \cite{Mao2015,Hul2018}   for each $\Delta$, we can  find a suitable  $l_{\Delta}$  that makes the $y_{\Delta, l_{\Delta}}(t)$ strongly converge to $y(t)$.  

Now let us describe  the construction of the truncated Euler-Maruyama method. First,  we assume  that there exists a function 
  $\psi: [1,\infty)\rightarrow(0,\infty) $, which   is strictly increasing and $\lim\limits_{R\rightarrow +\infty}\psi(R)=+\infty$,
 such that 
\begin{equation}\label{asc}
\sup\limits_{\mid y\mid\leq R}(\mid f(y)\mid \vee \mid g(y)\mid) \leq \psi (R), ~\forall R\geq 1\,. 
\end{equation}
We also assume that the  inverse $\psi^{-1}:[\psi(1),\infty)\rightarrow (0,\infty)$ of $\psi$ is also increasing.  At the same time, we define another strictly decreasing function $h:(0,1]\rightarrow[\psi(1),\infty)$, which satisfies
\begin{equation*}
\lim\limits_{\Delta\rightarrow0}h(\Delta)=+\infty, ~~~~~~~\Delta^{\frac{1}{4}}h(\Delta)\leq c_{0},
\end{equation*}
where $c_{0}$ is a positive constant which satisfies $c_{0}\geq1\vee \psi(1)$. Taking $l_{\Delta}=\psi^{-1}(h(\Delta))$, we get
\begin{equation}\label{eq5''}
\mathrm{d}y_{l_{\Delta}}(t)=f_{\Delta}(y_{l_{\Delta}}(t))\mathrm{d}t+g_{\Delta}(y_{l_{\Delta}}(t))\mathrm{d}W(t),
\end{equation}
where
\begin{equation*}
f_{\Delta}(x)=f\big((x\wedge\psi^{-1}(h(\Delta)))\frac{x}{\mid x\mid}\big),\quad    g_{\Delta}(x)=g\big((x\wedge\psi^{-1}(h(\Delta)))\frac{x}{\mid x\mid}\big).
\end{equation*}
The above equations together with \eqref{asc} imply that
\begin{equation}\label{boud}
\mid f_{\Delta}(x)\mid\leq h(\Delta),\quad \mid g_{\Delta} (x)\mid\leq h(\Delta).
\end{equation}
Then, applying the Euler-Maruyama method to \eqref{eq5''} yields
\begin{equation}\label{tre}
y_{\Delta}(t_{k+1})=y_{\Delta}(t_{k})+\Delta f_{\Delta}(y_{\Delta}(t_{k}))+g_{\Delta}(y_{\Delta}(t_{k}))\delta W_{k}, \quad k=0,1\cdots,N-1.
\end{equation}
Here, we set $y_{\Delta}(0)=y_{0}$.  
Define the continuous form of \eqref{tre} as
\begin{equation}\label{trec}
y_{\Delta}(t)=y_{0}+\int_{0}^{t}f_{\Delta}(\bar{y}_{\Delta}(s))\mathrm{d}s+\int_{0}^{t}g_{\Delta}(\bar{y}_{\Delta}(s))\mathrm{d}W(s),
\end{equation}
where $\bar{y}_{\Delta}(t)=y_{\Delta}(t_{k})$ for $t\in[t_{k},t_{k+1})$.

Next, we will prove the exponential integrability of both the exact solution $y(t)$ of SDE \eqref{eq5} and the numerical solution $y_{\Delta}(t)$. From now on, we always assume that $C$ stands generic  constants, whose values might change every time as it appears.
 
\begin{Theo}
Let  Assumption \ref{as1} hold  and let  $m<2$.  Then for any positive constant $q$, the exact solution $y(t)$ of \eqref{eq5} satisfies
\begin{equation}
\sup\limits_{t\in[0,T]}\EE( e^{q\mid y(t)\mid})<\infty.
\end{equation}
\end{Theo}
\textbf{Proof.}
An application of the  It\^{o} formula to  $y^{p}(t)$ yields 
\begin{equation*}
\begin{aligned}
 y^{p}(t)=y^{p}_{0}+\int_{0}^{t}py^{p-2}(s)[y(s)f(y(s))+\frac{p-1}{2}g^{2}(y(s))]\mathrm{d}s
 +\int_{0}^{t}py^{p-1}(s)g(y(s))\mathrm{d}W(s).
\end{aligned}
\end{equation*}
Here we only consider the case where p is a positive integer. Let $\tau_{R}=\inf\{t\in[0,T]: \mid y(t)\mid\geq R\}$. Then  we   derive that
\begin{equation*}
\begin{aligned}
\EE \big(y(t\wedge\tau_{R})^{2p}\big)&=y^{2p}_{0}+\EE \int_{0}^{t\wedge\tau_{R}}2py(s)^{2p-2}[y(s)f(y(s))+\frac{2p-1}{2}g^{2}(y(s))]\mathrm{d}s\\
&\leq y^{2p}_{0}+2Kp\EE \int_{0}^{t\wedge\tau_{R}}(2p)^{m}y(s)^{2p-2}\mathrm{d}s\\
&\leq y^{2p}_{0}+2 p\int_{0}^{t\wedge\tau_{R}}[\frac{1}{p}(2p)^{mp}+\frac{2p-2}{2p}\EE\big( y(s)^{2p}\big)]\mathrm{d}s\\
&\leq y^{2p}_{0}+2t(2p)^{mp}+Cp\int_{0}^{t}\EE \big(y(s\wedge\tau_{R})^{2p}\big)\mathrm{ds}.
\end{aligned}
\end{equation*}
Using the Gronwall inequality, we obtain
\begin{equation*}
\EE  \big(y(t\wedge\tau_{R})^{2p}\big)\leq\big(y_{0}^{2p}+2Kt(2p)^{mp}\big)e^{Cpt}\leq C_{t,y_{0}}^{p}(2p)^{mp}.
\end{equation*}
Besides, using the definition of $\tau_{R}$, we can infer
\begin{equation*}
\begin{aligned}
R^{2p}P(\tau_{R}\leq t)=\EE \big(y(\tau_{R})^{2p}I_{\{ \tau_{R}\leq t\}}\big)\leq \EE \big(  y(t\wedge\tau_{R})^{2p}\big)\leq C_{t,y_{0}}^{p}(2p)^{mp}.
\end{aligned}
\end{equation*} 
As a result, it can be derived  that 
\begin{equation*}
P(\tau_{\infty}>t)=1,
\end{equation*}
where $\tau_{\infty}=\lim\limits_{R\rightarrow \infty}\tau_{R}$. By means of the Fatou lemma, we get
\begin{equation*}
\EE  \big(y(t)^{2p}\big)\leq \varliminf\limits_{R\rightarrow \infty}\EE  \big(y(t\wedge\tau_{R})^{2p}\big)\leq C_{T,y(0)}^{p}(2p)^{mp}.
\end{equation*}
According to the Jensen inequality, we have
\begin{equation*}
\EE \mid y(t)\mid^{p}\leq \big(\EE (y(t)^{2p})\big)^{\frac{1}{2}}\leq C_{T,y_{0}}^{p}(2p)^{\frac{mp}{2}}\leq C_{T,y_{0}}^{p}p^{\frac{mp}{2}}.
\end{equation*}
It is a known fact that for any $q>0$, we have
\begin{equation*}
\EE e^{q\mid y(t)\mid}=\EE \sum\limits_{p=0}^{\infty}\frac{q^{p}\mid y(t)\mid^{p}}{p!}=\sum\limits_{p=0}^{\infty}\frac{q^{p}\EE \mid y(t)\mid^{p}}{p!}.
\end{equation*}
Finally, if $m<2$, by means of the Stirling formula, we can derive
\begin{equation*}
\EE ( e^{q\mid y(t)\mid})<\infty.
\end{equation*} 
This completes the proof of the theorem. 
\qed

\begin{Rem}
There are some papers concerned exponential integrability of exact solution. For example,  different kinds of exponential integrability are discussed in \cite{Hu2019, Hu1996} for SDE whose diffusion coefficient is a constant or grows linearly at most. The exponential integrability discussed in this paper is applicable to SDEs with exponentially growing drift and diffusion coefficients.
\end{Rem}

Next, we  prove that  the truncated Euler-Maruyama solutions are also  exponentially  integrable. But first, we need to state some lemmas.

\begin{Lem}
Let Assumption \ref{as1} hold. Then for arbitrary $\Delta\in(0,1]$, $p>2$ and $y\in\mathbb{R}$, we have
\begin{equation}
y f_{\Delta}(y)+\frac{p-1}{2}g_{\Delta}^{2}(y)\leq K(p^{m}+\mid y\mid).
\end{equation}
\end{Lem}
\textbf{Proof.}
When $\mid y\mid\leq\psi^{-1}(h(\Delta))$, $f_{\Delta}(y)=f(y),g_{\Delta}(y)=g(y)$, the conclusion is clear. When $\mid y\mid>\psi^{-1}(h(\Delta))$, 
we have
\begin{equation*}
\begin{aligned}
&yf_{\Delta}(y)+\frac{p-1}{2}g_{\Delta}^{2}(y)\\
\leq&\psi^{-1}(h(\Delta))\frac{y}{\mid y\mid}f_{\Delta}\big(\psi^{-1}(h(\Delta))\frac{y}{\mid y\mid}\big)+\frac{p-1}{2}g_{\Delta}^{2}\big(\psi^{-1}(h(\Delta))\frac{y}{\mid y\mid}\big)\\
&+\big(\frac{\mid y\mid}{\psi^{-1}(h(\Delta))}-1\big)\psi^{-1}(h(\Delta))\frac{y}{\mid y\mid}f_{\Delta}\big(\psi^{-1}(h(\Delta))\frac{y}{\mid y\mid}\big)\,. 
\end{aligned}
\end{equation*}
Since $yf(y)$ does not depend on the constant $p$, from condition \eqref{moas}, we can derive that there exists a positive constant $K$, independent of $p$ and $y$, such that $yf(y)\leq K$. Besides, we have $\psi^{-1}(h(\Delta))\geq\psi^{-1}(h(1))$. Then   the assertion is proved. 
\qed

\begin{Lem}
Let $Z\sim N(0,\sqrt{\Delta})$ be a one dimensional normal random variable. Then for arbitrary positive constant $\beta$, it holds 
\begin{equation}
\EE (e^{\beta \mid Z\mid})\leq 2e^{\frac{\beta^{2}\Delta}{2}}.
\end{equation}
\end{Lem}
\textbf{Proof.}
Since $Z$ is a normal random variable , it holds that
\begin{equation*}
\begin{aligned}
\EE (e^{\beta\mid Z\mid})&=\int_{\mathbb{R}}e^{\beta\mid x\mid}\frac{1}{\sqrt{2\pi\Delta}}e^{-\frac{x^{2}}{2\Delta}}\mathrm{d}x\\
&=\frac{2}{\sqrt{2\pi\Delta}}\int_{0}^{\infty}e^{\beta x}e^{-\frac{x^{2}}{2\Delta}}\mathrm{d}x\\
&=\frac{2}{\sqrt{2\pi}}e^{\frac{\beta^{2}\Delta}{2}}[\int_{-\beta\sqrt{\Delta}}^{0}e^{-\frac{x^{2}}{2}}\mathrm{d}x+\int_{0}^{\infty}e^{-\frac{x^{2}}{2}}\mathrm{d}x]\\
&=\frac{2}{\sqrt{2\pi}}e^{\frac{\beta^{2}\Delta}{2}}[\int_{0}^{\beta\sqrt{\Delta}}e^{-\frac{x^{2}}{2}}\mathrm{d}x+\int_{0}^{\infty}e^{-\frac{x^{2}}{2}}\mathrm{d}x]\\
&\leq\frac{2}{\sqrt{2\pi}}e^{\frac{\beta^{2}\Delta}{2}}[\bigg ( \int\!\!\!\!\int_{D} e^{-\frac{x^{2}+y^{2}}{2}}\mathrm{d}x\mathrm{d}y\bigg)^{\frac{1}{2}}+\int_{0}^{\infty}e^{-\frac{x^{2}}{2}}\mathrm{d}x]\\
&=\frac{2}{\sqrt{2\pi}}e^{\frac{\beta^{2}\Delta}{2}}\frac{\sqrt{\pi}}{\sqrt{2}}(\sqrt{1-e^{-\beta^{2}\Delta}}+1)\\
&\leq 2e^{\frac{\beta^{2}\Delta}{2}},
\end{aligned}
\end{equation*}
where $D=\{x\in\mathbb{R_{+}},y\in\mathbb{R_{+}}:x^{2}+y^{2}\leq 2\beta^{2}\Delta \}$.
\qed
\begin{Lem}
Under Assumption \ref{as1}, for any fixed step size $\Delta\in(0,1]$, the truncated Euler-Maruyama approximations for SDE \eqref{eq5} are exponentially  integrable, i.e.
\begin{equation*}
\EE(e^{q\mid y_{\Delta}(t)\mid})< \infty, ~~~\forall q>0.
\end{equation*}
\end{Lem}
\textbf{Proof.}
This can be easily proved by  combining  \eqref{boud} with  the above lemma.
\qed

\begin{Lem}\label{4.7}
The truncated Euler-Maruyama approximations defined by \eqref{trec} 
satisfy
\begin{equation}\label{3.11}
\EE  \mid y_{\Delta}(t)-\bar{y}_{\Delta}(t)\mid^{p}\leq C^{p}p^{p}\Delta^{\frac{p}{2}}(h(\Delta))^{p},~~\forall t\in[0,T], p>2.
\end{equation}
\end{Lem}
\textbf{Proof.}
 According to \eqref{trec}, for arbitrary $t\in[t_{k},t_{k+1})$, we have 
\begin{equation}\label{4aa}
\begin{aligned}
\EE\mid y_{\Delta}(t)-\bar{y}_{\Delta}(t)\mid^{p}=&\EE\mid \int_{t_{k}}^{t} f_{\Delta}(y_{\Delta}(t_{k}))\mathrm{d}s+\int_{t_{k}}^{t}g_{\Delta}(y_{\Delta}(t_{k}))\mathrm{d}W(s)\mid^{p}\\
\leq &2^{p-1}\big[\EE|\int_{t_{k}}^{t}f_{\Delta}(y_{\Delta}(t_{k}))\mathrm{d}s|^{p}+ \EE|\int_{t_{k}}^{t}g_{\Delta}(y_{\Delta}(t_{k}))\mathrm{d}W(s)|^{p}\big]\\
\leq&2^{p-1}\big[\Delta^{p-1}\int_{t_{k}}^{t}\EE\mid f_{\Delta}(y_{\Delta}(t_{k}))\mid^{p}\mathrm{d}s\\
&+\big( \frac{p(p-1)}{2}\big)^{\frac{p}{2}}\Delta^{\frac{p}{2}-1}\int_{t_{k}}^{t}\EE
\mid g_{\Delta}(y_{\Delta}(t_{k}))\mid^{p}\mathrm{d}s\big]\\
\leq&2^{p-1}\big[\Delta^{p}(h(\Delta))^{p}+\big( \frac{p(p-1)}{2}\big)^{\frac{p}{2}}\Delta^{\frac{p}{2}}(h(\Delta))^{p}\big]\\
\leq&C^{p}p^{p}\Delta^{\frac{p}{2}}(h(\Delta))^{p}.
\end{aligned}
\end{equation}
In the above inequalities, we have used the basic inequality $(a+b)^{p}\leq 2^{p-1}(a^{p}+b^{p}),\ \forall  a,b\geq0$, the inequality \eqref{boud},  and the  Burkholder-Davis-Gundy inequality  of It\^{o} stochastic integral.  
\qed

\begin{Rem}
Since $\Delta^{\frac{1}{4}}h(\Delta)\leq c_{0}$,  we can furthermore obtain
\begin{equation}
\EE\mid y_{\Delta}(t)-\bar{y}_{\Delta}(t)\mid^{p}\leq C^{p}p^{p}\Delta^{\frac{p}{4}},~~\forall t\in[0,T], p>2.
\end{equation}
\end{Rem}

\begin{Theo}
Let Assumption \ref{as1} hold, for any $q>0$, $m<2$, the truncated Euler-Maruyama approximations $y_{\Delta}(t)$ for SDE \eqref{eq5} satisfy
\begin{equation*}
\sup\limits_{\Delta\in(0,1]}\sup\limits_{t\in[0,T]}\EE(e^{q\mid y_{\Delta}(t)\mid})<\infty,~~~\forall q>0.
\end{equation*}
\end{Theo}
\textbf{Proof.}
By the It\^{o} formula and property of stochastic It\^{o} integral, we have
\begin{equation*}
\begin{aligned}
&\EE \big(y_{\Delta}(t)^{2p}\big)\\=&y^{2p}_{0}+E\int_{0}^{t}2py_{\Delta}(s)^{2p-2}[y_{\Delta}(s)f_{\Delta}(\bar{y}_{\Delta}(s))+\frac{2p-1}{2}g_{\Delta}^{2}(\bar{y}_{\Delta}(s))]\mathrm{d}s\\
=&y^{2p}_{0}+\EE \int_{0}^{t}2py_{\Delta}(s)^{2p-2}[\bar{y}_{\Delta}(s)f_{\Delta}(\bar{y}_{\Delta}(s))+\frac{2p-1}{2}g_{\Delta}^{2}(\bar{y}_{\Delta}(s))]\mathrm{d}s\\
&+\EE\int_{0}^{t}2py_{\Delta}(s)^{2p-2}[(y_{\Delta}(s)-\bar{y}_{\Delta}(s))f_{\Delta}(\bar{y}_{\Delta}(s))]\mathrm{d}s\\
=&y^{p}_{0}+A+B.
\end{aligned}
\end{equation*}
By means of the Young inequality $ab\leq \frac{a^{p}}{p}+\frac{b^{q}}{q}$, where $a,b>0, p>1, \frac{1}{p}+\frac{1}{q}=1$, we derive
\begin{equation*}
\begin{aligned}
A\leq& 2pK\EE \int_{0}^{t}[(2p)^{m}y_{\Delta}(s)^{2p-2}+y_{\Delta}(s)^{2p-2}\mid\bar{y}_{\Delta}(s)\mid]\mathrm{d}s\\
\leq&2Kp\big(\int_{0}^{t}\bigg[\frac{1}{p}(2p)^{mp}+\frac{2p-2}{2p}\EE \big(y_{\Delta}(s)^{2p}\big)+\frac{2p-2}{2p}\EE \big(y_{\Delta}(s)^{2p}\big)+\frac{1}{p}\big(\frac{1}{2}+\frac{1}{2}\EE (\bar{y}_{\Delta}(s)^{2p}) \big) \big]\mathrm{d}s\\
\leq&Ktp^{mp}+Cp\int_{0}^{t}\big[\EE \big(y_{\Delta}(s)^{2p}\big)+\EE \big(\bar{y}_{\Delta}(s)^{2p}\big)\big]\mathrm{d}s.
\end{aligned}
\end{equation*}
By using the Young inequality again,  we have
\begin{equation*}
\begin{aligned}
B\leq& (2p-2)\int_{0}^{t} \EE \big(y_{\Delta}(s)^{2p}\big)\mathrm{d}s+2\EE \int_{0}^{t}\mid y_{\Delta}(s)-\bar{y}_{\Delta}(s)\mid^{p}\mid f_{\Delta}(\bar{y}_{\Delta}(s))\mid^{p}\mathrm{d}s\\
\leq& Cp\int_{0}^{t} \EE \big(y_{\Delta}(s)^{2p}\big)\mathrm{d}s+2h(\Delta)^{p}\int_{0}^{t}\EE \mid y_{\Delta}(s)-\bar{y}_{\Delta}(s)\mid^{p}\mathrm{d}s\\
\leq&C p\int_{0}^{t} \EE \big(y_{\Delta}(s)^{2p}\big)\mathrm{d}s+2C^{p}p^{p}t(h(\Delta))^{2p}\Delta^{\frac{p}{2}}\\
\leq&Cp\int_{0}^{t} \EE \big(y_{\Delta}(s)^{2p}\big)\mathrm{d}s+C^{p}tp^{p}.
\end{aligned}
\end{equation*}
Consequently,  we get 
\begin{equation*}
\sup\limits_{s\in[0,t]}\EE\big( y_{\Delta}(s)^{2p}\big)\leq y^{2p}_{0}+C^{p}tp^{mp}+Cp\int_{0}^{t}\sup\limits_{u\in[0,s]}\EE  \big(y_{\Delta}(u)^{2p}\big)\mathrm{d}s.
\end{equation*}
Using the Gronwall inequality we have
\begin{equation*}
\sup\limits_{s\in[0,t]}\EE \big(y_{\Delta}(s)^{2p}\big)\leq (y^{2p}_{0}+C^{p}tp^{mp})e^{Cpt}.
\end{equation*}
Hence, we can get that for arbitrary $t>0$,
\begin{equation*}
\EE \big(y_{\Delta}(t)^{2p}\big)\leq C_{t,y_{0}}^{p}(p)^{mp}.
\end{equation*}
The Jensen inequality yields
\begin{equation*}
\EE \mid y_{\Delta}(t)\mid^{p}\leq C_{t,y_{0}}^{p}(p)^{\frac{mp}{2}}, m<2.
\end{equation*}
Now,  we can get the exponential integrability of truncated Euler-Maruyama approximation  for any $\Delta\in(0,1]$ by the Stirling formula.
\qed

Define
\begin{equation}
\tilde{\tau}_{R}=\inf\{ t>0: \mid y_{\Delta}(t)\mid\geq R\}\,.
\end{equation}
Since we already have the conclusion that both the exact solution and the numerical approximation are exponentially  integrable, 
  we can get
\begin{equation}
P(\tau_{R}\leq t )\leq \frac{C}{e^{qR}},~~ P(\tilde{\tau}_{R}\leq t )\leq \frac{C}{e^{qR}}.
\end{equation}
Thus under Assumption \eqref{moas},  the strong convergence of truncated Euler-Maruyama method is implied by   \cite[Theorem 3.5]{Mao2015}.
\begin{Lem}(\cite{Mao2015})
Let the coefficients of SDE \eqref{eq5} satisfy
\begin{equation*}
yf(y)+\frac{p-1}{2}g^{2}(y)\leq K(1+y^{2}),~\forall p>2\,. 
\end{equation*}
Then the truncated Euler-Maruyama approximation $y_{\Delta}(t)$ will strongly converge to the exact solution $y(t)$ of SDE \eqref{eq5}, i.e.
\begin{equation}
\lim\limits_{\Delta\rightarrow 0}\EE \mid y(t)-y_{\Delta}(t)\mid^{q}=0,  ~\forall q>0.
\end{equation}
\end{Lem}

\begin{Theo}\label{3a}
Let Assumption \ref{as1} hold. Then the numerical solution $x_{\Delta}(t)=e^{y_{\Delta}(t)}$ will strongly converge to the exact solution $x(t)$ of SDE \eqref{eq1}, i.e.
\begin{equation}
\lim\limits_{\Delta\rightarrow 0}\EE \mid x(t)-x_{\Delta}(t)\mid^{q}=0, ~\forall q>0.
\end{equation}
More importantly, the numerical solution $x_{\Delta}(t)$  maintains the positiveness of the original solution.
\end{Theo}
\textbf{Proof.}
By means of the mean value theorem, we obtain
\begin{equation*}
\begin{aligned}
\EE \mid x(t)-x_{\Delta}(t)\mid^{q}=&\EE \mid e^{y(t)}-e^{y_{\Delta}(t)}\mid^{q}\\
\leq &\EE \mid e^{y(t)}+e^{y_{\Delta}(t)}\mid^{q}\mid y(t)-y_{\Delta}(t)\mid^{q}\\
\leq &(\EE \mid e^{y(t)}+e^{y_{\Delta}(t)}\mid^{2q})^{\frac{1}{2}}(\EE \mid y(t)-y_{\Delta}(t)\mid^{2q})^{\frac{1}{2}}.
\end{aligned}
\end{equation*}
The convergence result  can be derived by using exponential integrability results and the above lemma.
\qed

\begin{Exe}\label{a}
Consider the CEV process
\begin{equation}\label{cev}
\mathrm{d}x(t)=\kappa(\lambda-x(t))\mathrm{d}t+\theta x(t)^{\alpha}\mathrm{d}W(t),
\end{equation}
where constants $\kappa,\lambda,\theta>0$, and $\alpha\in(\frac{1}{2},1)$. 
We already know that this SDE has a positive solution, and the boundary  $0$ is always unattainable (see \cite{Andersen2007}). Applying the It\^{o} formula  to  $y(t)=\ln x(t)$, we have
\begin{equation}\label{tcev}
\mathrm{d}y(t)=[\kappa\lambda e^{-y(t)}-\frac{\theta^{2}}{2}e^{-2(1-\alpha)y(t)}-\kappa]\mathrm{d}t+\theta e^{-(1-\alpha)y(t)}\mathrm{d}W(t).
\end{equation} 
Denote
\begin{equation*}
A:= y[\kappa\lambda e^{-y}-\frac{\theta^{2}}{2}e^{-2(1-\alpha)y}-\kappa]+\frac{p-1}{2}\theta^{2}e^{-2(1-\alpha)y}.
\end{equation*}
Since $2(1-\alpha)<1$, we have
\begin{equation*}
A\leq Kp 
\end{equation*}
for $y>0$.  Since
\begin{equation*}
\lim\limits_{y\rightarrow -\infty}\frac{\frac{-\theta^{2}}{2}ye^{-2(1-\alpha)y}}{\frac{\kappa\lambda}{2}e^{-y}}=0,
\end{equation*}
there exists a constant $-M$, such that 
\begin{equation*}
\begin{aligned}
A\leq -\kappa\lambda e^{-y}+\frac{\kappa\lambda}{2}e^{-y}+\frac{\kappa}{2(1-\alpha)}e^{-2(1-\alpha)y}+\frac{p-1}{2}\theta^{2}e^{-2(1-\alpha)y} ,
\end{aligned}
\end{equation*}
for $y<-M\wedge -1$. Denote
\begin{equation*}
H(y)= -\kappa\lambda e^{-y}+\frac{\kappa\lambda}{2}e^{-y}+\frac{\kappa}{2(1-\alpha)}e^{-2(1-\alpha)y}+\frac{p-1}{2}\theta^{2}e^{-2(1-\alpha)y}.
\end{equation*}
Taking derivative of the above equality yields 
\begin{equation*}
\begin{aligned}
H'(y)=&\frac{\kappa\lambda}{2}e^{-y}-\kappa e^{-2(1-\alpha)y}-(p-1)(1-\alpha)\theta^{2}e^{-2(1-\alpha)y}\\
=&e^{-y}\bigg(\frac{\kappa\lambda}{2}-\kappa e^{[1-2(1-\alpha)]y}-(p-1)(1-\alpha)\theta^{2}e^{[1-2(1-\alpha)]y}\bigg)\,. 
\end{aligned}
\end{equation*}
Solving $H'(y_{0})=0$ we  obtain  
\begin{equation*}
y_{0}=\frac{1}{2\alpha-1}\ln\frac{\kappa\lambda}{2\kappa+2(p-1)(1-\alpha)\theta^{2}}.
\end{equation*}
Therefore there exists a positive constants $K$ such that
\begin{equation*}
A\leq Kp^{1+\frac{2(1-\alpha)}{2\alpha-1}}.
\end{equation*}
 If $\alpha \in(\frac{3}{4},1)$, the coefficients of \eqref{tcev} satisfies  condition \eqref{moas} with $m =1+\frac{2(1-\alpha)}{2\alpha-1}<2$.
This means that the  coefficients of the transformed SDE satisfy Assumption \ref{as1}.
We can apply then Theorem \ref{3a} to the CEV process.
\end{Exe}

\begin{Exe}\label{b}
Consider the generalized Ait-Sahalia model (discussed in \cite{szpruch2011})
\begin{equation}\label{air}
\mathrm{d}x(t)=[a_{-1}x(t)^{-1}-a_{0}+a_{1}x(t)-a_{2}x(t)^{r}]\mathrm{d}t+\sigma x(t)^{\rho}\mathrm{d}W(t),
\end{equation}
where the constants $a_{-1},a_{0},a_{1},a_{2},\sigma>0$, and $r>1,\rho>1$. Since the coefficients satisfy the local Lipschitz condition, using the standard analytical technique, we know that there is  a maximal local solution $x(t)$ for $t\in[0,\tau_{s})$, where $\tau_{s}$ is the stopping time when the solution explodes or first reaches boundary zero.   
Since 
\begin{equation*}
\begin{aligned}
p(x)=&\int_{c}^{x}\exp\{-2\int_{c}^{\xi}\frac{a_{-1}u^{-1}-a_{0}+a_{1}u-a_{2}u^{r}}{\sigma^{2}u^{2\rho}}\mathrm{d}u\}\mathrm{d}\xi,
\end{aligned}
\end{equation*} 
here $c>0$, and
\begin{equation*}
\frac{a_{-1}u^{-1}-a_{0}+a_{1}u-a_{2}u^{r}}{\sigma^{2}u^{2\rho}}\sim_{0+}\frac{a_{-1}}{\sigma^{2}u^{2\rho+1}}.
\end{equation*}
However, when $x\rightarrow0+$, it can be derived that
\begin{equation*}
\begin{aligned}
\int_{c}^{x}\exp\{-2\int_{c}^{\xi}\frac{a_{-1}}{\sigma^{2}u^{2\rho+1}}\mathrm{d}u\}\mathrm{d}\xi=C\int_{c}^{x}e^{\frac{a_{-1}}{\rho\sigma^{2}}\xi^{-2\rho}}\mathrm{d}\xi\rightarrow-\infty.
\end{aligned}
\end{equation*}
So $p(0+)=-\infty$.  On the other hand,  
\begin{equation*}
\frac{a_{-1}u^{-1}-a_{0}+a_{1}u-a_{2}u^{r}}{\sigma^{2}u^{2\rho}}\sim_{\infty}\frac{-a_{2}u^{r}}{\sigma^{2}u^{2\rho}}.
\end{equation*}
When $x\rightarrow+\infty$, we get
\begin{equation*}
\begin{aligned}
\int_{c}^{x}\exp\{-2\int_{c}^{\xi}\frac{-a_{2}u^{r}}{\sigma^{2}u^{2\rho}}\mathrm{d}u\}\mathrm{d}\xi=C\int_{c}^{x}e^{\frac{2a_{2}}{(r-2\rho+1)\sigma^{2}}\xi^{r-2\rho+1}}\mathrm{d}\xi\rightarrow+\infty.
\end{aligned}
\end{equation*}
Hence $p(+\infty)=+\infty$.  We also know that $p(0+)=-\infty\Rightarrow v(0+)=+\infty, p(+\infty)=+\infty\Rightarrow v(+\infty)=+\infty$. An application of the Feller criteria yields  
\[
P\big(\lim\limits_{k\rightarrow \infty}\inf\{t\geq 0: \frac{1}{k}<x(t)<k\}=\infty\big)=1\,,
\]
which means this SDE has a  unique positive global solution and both the boundaries are unattainable.
The application of  the It\^{o} formula with $y(t)=\ln x(t)$ yields  
\begin{equation}\label{tair}
\begin{aligned}
\mathrm{d}y(t)=f(y(t))\mathrm{d}t+g(y(t))\mathrm{d}W(t),
\end{aligned}
\end{equation}
where 
\[
f(y)=a_{-1}e^{-2y}-a_{0}e^{-y}+a_{1}-a_{2}e^{(r-1)y}-\frac{\sigma^{2}}{2}e^{2(\rho-1)y}, \quad g(y)=\sigma e^{(\rho-1)y}\,.
\]
Denote 
\begin{equation*}
B:=yf(y)+\frac{p-1}{2}g^{2}(y).
\end{equation*}
When $r>2\rho-1$, we have $B\leq Kp$ for $y<0$, and  
\begin{equation*}
B\leq C+\frac{a_{1}}{2(\rho-1)}e^{2(\rho-1)y}-a_{2}e^{(r-1)y}+\frac{p-1}{2}\sigma^{2}e^{2(\rho-1)y},
\end{equation*}
when $y>1$. Let 
\begin{equation*}
 H_{1}(y)=C+\frac{a_{1}}{2(\rho-1)}e^{2(\rho-1)y}-a_{2}e^{(r-1)y}+\frac{p-1}{2}\sigma^{2}e^{2(\rho-1)y}.
\end{equation*}
We have 
\begin{equation*}
\begin{aligned}
H_{1}'(y)=&a_{1}e^{2(\rho-1)y}-(r-1)a_{2}e^{(r-1)y}+\frac{p-1}{2}\sigma^{2}e^{2(\rho-1)y}\\
=&e^{(r-1)y}\big( a_{1}e^{-(r-2\rho+1)y}-(r-1)a_{2}+\frac{p-1}{2}\sigma^{2}e^{-(r-2\rho+1)y}\big).
\end{aligned}
\end{equation*}
The above derivative satisfies $H_{1}'(y_{1})=0$,  when  
\begin{equation*}
y_{1}=\frac{1}{2-2\rho+1}\ln\frac{\sigma^{2}(p-1)(\rho-1)+a_{1}}{a_{2}(r-1)}.
\end{equation*}
Therefore we can get that
\begin{equation*}
B\leq K p^{1+\frac{2(\rho-1)}{r-2\rho-1+1}}.
\end{equation*}
When $r>4\rho-3$, the coefficients of SDE \eqref{tair} satisfy condition \eqref{moas} with $m=1+\frac{2(\rho-1)}{r-2\rho-1+1}<2$.  Thus,  we can find a positive constant $K$  such that Assumption \ref{as1} holds. This means that  Theorem \ref{3a} can be applied to the generalized Ait-Sahalia model \eqref{air}.
\end{Exe}

\begin{Rem}
For the CEV process \eqref{cev} and the generalized Ait-Sahalia model \eqref{air}, we know that $\EE  \big(x(t)^{q}\big)<\infty$, $\EE\big( x(t)^{-q}\big)<\infty$, for $\forall q>0$, with $\alpha\in(\frac{1}{2},1)$ and $r>2\rho-1$  respectively. Since   $x(t)=e^{y(t)}$,   this is equivalent to  the conclusion $\EE  e^{q\mid y(t)\mid}<\infty$, $\forall q>0$.  
In the above arguments,   the exponential integrability for the equations \eqref{tcev} and \eqref{tair} requires that parameters in coefficients of those SDEs satisfy $\alpha\in(\frac{3}{4},1)$ and $r>4\rho-3$, respectively. 
In fact, to prove Theorem \ref{3a}, we only  need to prove $\EE e^{qy(t)}<\infty$, $\forall q>0$. It is interesting to know if the condition \eqref{moas} can be weakened for the  convergence of our numerical method.
\end{Rem}
Moreover, for the CIR process \eqref{ciro}, after using logarithmic transformation, we obtain  SDE \eqref{cirt}. The coefficients of \eqref{cirt} does not satisfy the condition \eqref{moas}. This can be easily seen from the following computations. 
Denote
\begin{equation*}
A_{0}:=yf(y)+\frac{p-1}{2}g^{2}(y)=(\kappa\lambda-\frac{\theta^{2}}{2})ye^{-y}-\kappa y+\frac{p-1}{2}\theta^{2}e^{-y}.
\end{equation*}
 When $y>0$, it follows that $A_{0}\leq Kp$. But when $y<0$, denote
 \begin{equation*}
 H(y)=(\kappa\lambda-\frac{\theta^{2}}{2})ye^{-y}-\kappa y+\frac{p-1}{2}\theta^{2}e^{-y}.
 \end{equation*}
 The derivative has the following form
 \begin{equation*}
 \begin{aligned}
 H'(y)=&(\kappa\lambda-\frac{p}{2}\theta^{2})e^{-y}-(\kappa\lambda-\frac{\theta^{2}}{2})ye^{-y}-\kappa \\
 =&e^{-y}[ \kappa\lambda-\frac{p}{2}\theta^{2}-(\kappa\lambda-\frac{\theta^{2}}{2})y-\kappa e^{y}].
 \end{aligned}
 \end{equation*}
 If $\kappa\lambda<\frac{p}{2}\theta^{2}$,    then we have 
 \begin{equation*}
 \kappa\lambda-\frac{p}{2}\theta^{2}-(\kappa\lambda-\frac{\theta^{2}}{2})y\leq 0,
 \end{equation*}
 for $y\leq y_{1}$, $y_{1}=[\kappa\lambda-\frac{p}{2}\theta^{2}]/[\kappa\lambda-\frac{\theta^{2}}{2}]$.  It is easy to see  that $H'(y_{1})<0$. We assume $H'(y_{0})=0$. Then there exists a $y_{2}\in (y_{0},y_{1})$, such that
\begin{equation*}
H(y_{2})>H(y_{1})=(\kappa\lambda-\frac{\theta^{2}}{2})e^{c_{2}p-c_{1}}+\kappa(c_{2}p-c_{1})> C e^{c_{2}p},
\end{equation*}
where $c_{2}=\frac{\theta^{2}}{2\kappa\lambda-\theta^{2}}$ and $c_{1}=\frac{2\kappa\lambda}{2\kappa\lambda-\theta^{2}}$.  

To cover the case of CIR and other more general SDE, we introduce  the following assumption.
\begin{Assu}\label{as4.2}
Assume that the coefficients of  \eqref{eq5} satisfy
\begin{equation}\label{weakcd1}
yf(y)+\frac{p-1}{2}g^{2}(y)\leq K,\quad \forall y\in \RR\,, p>0\,, 
\end{equation}
where the positive constant $K=K_p$ depends only 
on $p$, but independent of
$y$. And  there exists a positive constant $K_{1}=K_{1,p}$ which also depends only on $p$, such that
\begin{equation}\label{weakcd2}
e^{2y}[f(y)+\frac{p}{2}g^{2}(y)]\leq K_{1} (1+e^{2y}),\quad   \forall y\in \RR\,, p>0.
\end{equation}
\end{Assu}
From  \cite{Mao2015}, we know that under condition \eqref{weakcd1}, the truncated Euler-Maruyama numerical solution $y_{\Delta}(t)$ strongly converges to the solution of SDE \eqref{eq5} and $P(\tilde{\tau}_{\infty}>T)=0$. 
Since
\begin{equation*}
e^{2y}[f(y)+\frac{p}{2}g^{2}(y)]=e^{y}b(e^{y})+\frac{p-1}{2}\sigma^{2}(e^{y}),
\end{equation*}
it is easy to see  that the condition \eqref{weakcd2} becomes 
\begin{equation*}
xb(x)+\frac{p-1}{2}\sigma^{2}(x)\leq K_{1}(1+x^{2})\,, 
\end{equation*}
which implies  $\EE\big( x(t)^{p})<\infty$, $\forall p>0$, and consequently, 
 $\EE  e^{py(t)}<\infty$. In order to get strong convergence result, we  need to prove $\EE  e^{py_{\Delta}(t)}<\infty$, and this is what the next theorem says.

\begin{Theo}
If the coefficients of SDE \eqref{eq5} satisfy Assumption \ref{as4.2}, then the truncated Euler-Maruyama approximation $y_{\Delta}(t)$ is exponentially  integrable, i.e.
\begin{equation*}
\sup\limits_{\Delta\in(0,1]}\sup\limits_{t\in[0,T]}\EE  e^{py_{\Delta}(t)}<\infty.
\end{equation*}
\end{Theo}
\textbf{Proof.}
By applying It\^{o} formula, we derive
\begin{equation*}
\begin{aligned}
\EE  e^{py_{\Delta}(t\wedge \tilde{\tau}_{R})}
=&e^{py_{0}}+\EE \int_{0}^{t\wedge\tilde{\tau}_{R}}pe^{(p-2)y_{\Delta}(s)}[e^{2y_{\Delta}(s)}\big( f_{\Delta}(\bar{y}_{\Delta}(s))+\frac{p}{2}g_{\Delta}^{2}(\bar{y}_{\Delta}(s))\big)]\mathrm{d}s\\
=&e^{py_{0}}+\EE \int_{0}^{t\wedge\tilde{\tau}_{R}}pe^{py_{\Delta}(s)-2\bar{y}_{\Delta}(s)}[e^{2\bar{y}_{\Delta}(s)}\big( f_{\Delta}(\bar{y}_{\Delta}(s))+\frac{p}{2}g_{\Delta}^{2}(\bar{y}_{\Delta}(s))\big)\mathrm{d}s\\
\leq&e^{py_{0}}+K_{1}p\EE \int_{0}^{t\wedge\tilde{\tau}_{R}}[e^{py_{\Delta}(s)}+e^{(p-2)y_{\Delta}(s)}e^{2\big(y_{\Delta}(s)-\bar{y}_{\Delta}(s)\big)}]\mathrm{d}s\\
\leq&e^{py_{0}}+Cp\int_{0}^{t\wedge\tilde{\tau}_{R}}\EE e^{py_{\Delta}(s)}\mathrm{d}s+2K_{1}\EE \int_{0}^{t\wedge\tilde{\tau}_{R}}e^{p\big(y_{\Delta}(s)-\bar{y}_{\Delta}(s)\big)}\mathrm{d}s.
\end{aligned}
\end{equation*}
For the last item on the right side of the above inequality, we have
\begin{equation*}
\begin{aligned}
\EE e^{p\big(y_{\Delta}(t)-\bar{y}_{\Delta}(t)\big)}\leq &\EE e^{p\mid y_{\Delta}(t)-\bar{y}_{\Delta}(t)\mid}\\
\leq&\EE e^{p[\int_{t_{k}}^{t}\mid f_{\Delta}(y_{\Delta}(t_{k}))\mid\mathrm{d}s+\mid g_{\Delta}(y_{\Delta}(t_{k}))\mid\cdot\mid \int_{t_{k}}^{t}\mathrm{d}W(s)\mid]}\\
\leq&\EE e^{ph(\Delta)\big( \Delta+\mid W(t)-W(t_{k})\mid\big)}\\
\leq&e^{Cp^{2} h^{2}(\Delta)\Delta}\\
\leq&e^{Cp^{2}\Delta^{\frac{1}{2}}},
\end{aligned}
\end{equation*}
for $t\in[t_{k},t_{k+1})$. Thus, according to the Gronwall inequality, we derive
\begin{equation*}
\EE e^{py_{\Delta}(t\wedge \tilde{\tau}_{R})}\leq(e^{py_{0}}+2K_{1}te^{Cp^{2}})e^{Cpt}.
\end{equation*}
Since $P(\tilde{\tau}_{\infty}>t)=0$ for any $t>0$, an application of the Fatou lemma yields
\begin{equation*}
\sup\limits_{t\in[0,T]}\EE e^{py_{\Delta}(t)}<\infty, ~\forall \Delta\in(0,1].
\end{equation*} 
This proves the theorem. 
\qed

\begin{Theo}\label{3b}
If the coefficients of SDE \eqref{eq5} satisfy Assumption \ref{as4.2}, then for any $q>0$, the logTE numerical approximation $x_{\Delta}(t)=e^{y_{\Delta}(t)}$ strongly converges to exact solution of SDE \eqref{eq1}, i.e.
\begin{equation*}
\lim\limits_{\Delta\rightarrow 0}\EE \mid x(t)-x_{\Delta}\mid^{q}=0.
\end{equation*}
\end{Theo}
\textbf{Proof.}
This can be easily proved by using similar techniques to the proof of Theorem \ref{3a}. 
\qed

\begin{Rem}
It is not difficult to see  that the coefficients of SDE \eqref{cirt} satisfy Assumption \ref{as4.2}, therefore, Theorem \ref{3b} can be applied to CIR process \eqref{ciro}.
\end{Rem}

\section{Convergence rate of numerical scheme}
Although we have obtained  the strong convergence theorem for the logTE scheme, we still do not know how fast the numerical solution converges. In this section  we  impose stronger  conditions  on the coefficients than those  in the previous section to get the convergence rate. 
\begin{Assu}\label{as2r}
Let the coefficients of SDE \eqref{eq5}  satisfy the following conditions: there exist positive constants $K=K_{m}$, $K_{1}=K_{1,q}$ and $K_{2}$,  where $K_{1}$ only depends on $q$, $K_{2}$ is independent of $y_{1},y_{2}$, and $K$ only depends on $m$, such that, for any $y, y_{1},y_{2}\in \mathbb{R}$, we have
\begin{equation}
(y_{1}-y_{2})\big(f(y_{1})-f(y_{2})\big)+\frac{q-1}{2}\mid g(y_{1})-g(y_{2})\mid^{2}\leq K_{1}\mid y_{1}-y_{2}\mid^{2},~\forall q>2,
\end{equation}
\begin{equation}
\mid f(y_{1})-f(y_{2})\mid\vee\mid g(y_{1})-g(y_{2})\mid\leq K_{2}(e^{\beta\mid y_{1}\mid}+e^{\beta\mid y_{2}\mid})\mid y_{1}-y_{2}\mid,
\end{equation}
\begin{equation}\label{moas1}
yf(y)+\frac{p-1}{2}g^{2}(y)\leq Kp^{m},~1\leq m<2, \forall p>2 .
\end{equation}
\end{Assu}

The main theorem of this section is the following convergence rate of the  numerical scheme constructed  in the previous section.
\begin{Theo}\label{4.1}
Let the coefficients of SDE \eqref{eq5} satisfy Assumption \ref{as2r}. If define $\psi(u)=He^{(\beta+1)u}, u\geq 1$, $h(\Delta)=H_{1}\Delta^{-\epsilon}$, $x_{\Delta}(t)=e^{y_{\Delta}(t)}$, then for any $q_{1}>0$ and  $\epsilon\in(0,\frac{1}{4})$, it holds that
\begin{equation}
\begin{cases}
\EE \mid y(T)-y_{\Delta}(T)\mid^{q_{1}}\leq C_{\vare, q_1} \Delta^{\frac{q_{1}(1-2\epsilon)}{2}}\,;   \\
\EE \mid x(T)-x_{\Delta}(T)\mid^{q_{1}}\leq  C_{\vare, q_1}\Delta^{\frac{q_{1}(1-2\epsilon)}{2}},\\
\end{cases} 
\end{equation}
where $y(t)$, $x(t)$ are exact solutions of SDE \eqref{eq5}  and SDE \eqref{eq1} respectively .
\end{Theo}
\textbf{Proof.}
Denote 
\[
\nu_{R}=\tau_{R}\wedge \tilde{\tau}_{R}\,,
\quad e_{\Delta}(t)=y(t)-y_{\Delta}(t)\,.
\]
 An application of the It\^{o} formula to $e_{\Delta}(t)^{2q_{1}}$ 
 yields 
\begin{equation*}
\begin{aligned}
&\EE e_{\Delta}(t\wedge\nu_{R})^{2q_{1}}\\
=&\EE  \int_{0}^{t\wedge\nu_{R}}2q_{1}e_{\Delta}(s)^{2q_{1}-2}\{e_{\Delta}(s)[f(y(s))-f_{\Delta}(\bar{y}_{\Delta}(s))]
+\frac{2q_{1}-1}{2}\mid g(y(s))-g_{\Delta}(\bar{y}_{\Delta}(s))\mid^{2}\}\mathrm{d}s\\
\leq&\EE  \int_{0}^{t\wedge\nu_{R}}2q_{1}e_{\Delta}(s)^{2q_{1}-2}\{e_{\Delta}(s)[f(y(s))-f(y_{\Delta}(s))]
+\frac{(2q-1)}{2}\mid g(y(s))-g(y_{\Delta}(s))\mid^{2}\}\mathrm{d}s\\
&+\EE  \int_{0}^{t\wedge\nu_{R}}2q_{1}e_{\Delta}(s)^{2q_{1}-2}\{e_{\Delta}(s)[f(y_{\Delta}(s))-f_{\Delta}(\bar{y}_{\Delta}(s))]\\
&+\frac{(2q_{1}-1)(2q-1)}{2(2q-2q_{1})}\mid g(y_{\Delta}(s))-g_{\Delta}(\bar{y}_{\Delta}(s))\mid^{2}\}\mathrm{d}s\\
=&H_{1}+H_{2},
\end{aligned}
\end{equation*}
where $ q_{1}$ is positive integer   satisfying  $q_{1}<q$. It follows from Assumption \ref{as2r} that
\begin{equation*}
H_{1}\leq C\int_{0}^{t\wedge\nu_{R}}\EE e_{\Delta}(s)^{2q_{1}}\mathrm{d}s.
\end{equation*}
Besides, we can derive that
\begin{equation*}
\begin{aligned}
H_{2}
\leq &\EE  \int_{0}^{t\wedge\nu_{R}}2q_{1}e_{\Delta}(s)^{2q_{1}-2}\bigg\{e_{\Delta}(s)[f(y_{\Delta}(s))-f_{\Delta}(y_{\Delta}(s))]\\
&+\frac{(2q-1)(2q_{1}-1)}{2q-2q_{1}}\mid g(y_{\Delta}(s))-g_{\Delta}(y_{\Delta}(s))\mid^{2}\bigg\}\mathrm{d}s\\
&+\EE  \int_{0}^{t\wedge\nu_{R}}2q_{1}e_{\Delta}(s)^{2q_{1}-2}\bigg\{e_{\Delta}(s)[f_{\Delta}(y_{\Delta}(s))-f_{\Delta}(\bar{y}_{\Delta}(s))]\\
&+\frac{(2q-1)(2q_{1}-1)}{2q-2q_{1}}\mid g_{\Delta}(y_{\Delta}(s))-g_{\Delta}(\bar{y}_{\Delta}(s))\mid^{2}\bigg\}\mathrm{d}s\\
=&H_{21}+H_{22}.
\end{aligned}
\end{equation*}
Using the Young inequality twice, we have 
\begin{equation*}
\begin{aligned}
H_{21}\leq &\EE \int_{0}^{t\wedge\nu_{R}}2q_{1}e_{\Delta}(s)^{2q_{1}-2}\bigg\{\frac{1}{2}e_{\Delta}(s)^{2}+\frac{1}{2}\mid f(y_{\Delta}(s))-f_{\Delta}(y_{\Delta}(s))\mid^{2} \\
&+\frac{(2q-1)(2q_{1}-1)}{2q-2q_{1}}\mid  g(y_{\Delta}(s))-g_{\Delta}(y_{\Delta}(s))\mid^{2}\bigg\}\\
\leq & Cq_{1}\int_{0}^{t\wedge\nu_{R}}\EE  e_{\Delta}(s)^{2q_{1}}\mathrm{d}s+C\EE \int_{0}^{t\wedge\nu_{R}}\big[\mid  f(y_{\Delta}(s))-f_{\Delta}(y_{\Delta}(s)) \mid^{2q_{1}}\\
&+\mid  g(y_{\Delta}(s))-g_{\Delta}(y_{\Delta}(s))\mid^{2q_{1}}\big]\mathrm{d}s\,. 
\end{aligned}
\end{equation*}
Denote $\pi(y)=\big(\mid\ y\mid\wedge\psi^{-1}(h(\Delta))\big)\frac{y}{\mid y\mid}$. Obviously,  we have $\mid\pi(y)\mid\leq\mid y\mid$.  We have  
\begin{equation*}
\begin{aligned}
&\EE \int_{0}^{t\wedge\nu_{R}}\big[\mid  f(y_{\Delta}(s))-f_{\Delta}(y_{\Delta}(s)) \mid^{2q_{1}}+\mid  g(y_{\Delta}(s))-g_{\Delta}(y_{\Delta}(s))\mid^{2q_{1}}\big]\mathrm{d}s\\
\leq & C\int_{0}^{t\wedge\nu_{R}}\EE \bigg((e^{2\beta q_{1}\mid y_{\Delta}(s)\mid}+e^{2\beta q_{1} \mid \pi(y_{\Delta}(s))\mid})\mid y_{\Delta}(s)-\pi(y_{\Delta}(s))\mid^{2q_{1}}\bigg)\mathrm{d}s\\
\leq &C\int_{0}^{t\wedge\nu_{R}}\bigg(\EE (e^{2\beta q_{1}\mid y_{\Delta}(s)\mid}+e^{2\beta q_{1} \mid \pi(y_{\Delta}(s))\mid})^{\varrho}\bigg)^{\frac{1}{\varrho}}\bigg(\EE \mid y_{\Delta}(s)\!-\!\pi(y_{\Delta}(s))\mid^{\frac{2\varrho q_{1}}{\varrho-1}}\bigg)^{\frac{\varrho-1}{\varrho}}\mathrm{d}s\\
\leq&C\int_{0}^{t\wedge\nu_{R}}\bigg(\EE \mid y_{\Delta}(s)-\pi(y_{\Delta}(s))\mid^{\frac{2\varrho q_{1}}{\varrho-1}}\bigg)^{\frac{\varrho-1}{\varrho}}\mathrm{d}s,
\end{aligned}
\end{equation*}
where $\varrho>1$. For $\varrho_{1}>1$, it is easy to derive 
\begin{equation*}
\begin{aligned}
&\bigg(\EE \mid y_{\Delta}(s)-\pi(y_{\Delta}(s))\mid^{\frac{2\varrho q_{1}}{\varrho-1}}\bigg)^{\frac{\varrho-1}{\varrho}}\\
\leq&C\bigg(\EE \big[I_{\{\mid y_{\Delta}(s)\mid>\psi^{-1}(h(\Delta))\}}\mid y_{\Delta}(s)\mid^{\frac{2\varrho q_{1}}{\varrho-1}}\big]\bigg)^{\frac{\varrho-1}{\varrho}}\\
\leq&C\bigg(P\big[\mid y_{\Delta}(s)\mid>\psi^{-1}(h(\Delta))\big]\bigg)^{\frac{(\varrho_{1}-1)(\varrho-1)}{\varrho\varrho_{1}}}\bigg(\EE \mid y_{\Delta}(s)\mid^{\frac{2\varrho \varrho_{1}q_{1}}{\varrho-1}}\bigg)^{\frac{\varrho-1}{\varrho\varrho_{1}}}\\
\leq&C\bigg(\frac{\EE e^{p\mid y_{\Delta}(s)\mid}}{e^{p\psi^{-1}(h(\Delta))}}\bigg)^{\frac{(\varrho_{1}-1)(\varrho-1)}{\varrho\varrho_{1}}}\\
\leq&Ce^{-\frac{p(\varrho_{1}-1)(\varrho-1)}{\varrho\varrho_{1}}\psi^{-1}(h(\Delta))}.
\end{aligned}
\end{equation*}
On the other hand,    using Lemma \ref{4.7} and then the Young inequality, we get   
\begin{equation*}
H_{22}\leq C\int_{0}^{t\wedge\nu_{R}}\EE e_{\Delta}(s)^{2q_{1}}\mathrm{d}s+C\Delta^{q_{1}}\big(h(\Delta)\big)^{2q_{1}}.
\end{equation*}
Applications of the Gronwall inequality and the Jensen inequality yields 
\begin{equation*}
\EE \mid e_{\Delta}(t\wedge\nu_{R})\mid^{q_{1}} \leq C(\Delta^{\frac{\epsilon p(\varrho_{1}-1)(\varrho-1)}{2\varrho\varrho_{1}(\beta+1)}}+\Delta^{\frac{q_{1}(1-2\epsilon)}{2}}).
\end{equation*}
Since $P(\nu_{\infty}>t)=1$, we can apply  the Fatou lemma to  get
\begin{equation*}
\EE \mid e_{\Delta}(t)\mid^{q_{1}}\leq C(\Delta^{\frac{\epsilon p(\varrho_{1}-1)(\varrho-1)}{2\varrho\varrho_{1}(\beta+1)}}+\Delta^{\frac{q_{1}(1-2\epsilon)}{2}}).
\end{equation*}
Taking $p$ sufficiently large such that $\frac{ \epsilon p(\varrho_{1}-1)(\varrho-1)}{2\varrho\varrho_{1}(\beta+1)}\geq\frac{q_{1}(1-2\epsilon)}{2}$  we have  
\begin{equation*}
\EE \mid e_{\Delta}(t)\mid^{q_{1}}\leq C\Delta^{\frac{q_{1}(1-2\epsilon)}{2}}, \epsilon\in(0,\frac{1}{4}),\forall q_{1}>2.
\end{equation*}
Finally,   using the mean value theorem and the fact that both the exact solution and the numerical solution are exponentially  integrable, we derive
\begin{equation*}
\EE \mid x(t)-x_{\Delta}(t)\mid^{q_{1}} \leq C\Delta^{\frac{q_{1}(1-2\epsilon)}{2}}, \epsilon\in(0,\frac{1}{4}),\forall q_{1}>2.
\end{equation*}
Although in the beginning  of the proof, we assume  $q_{1}>2$, other cases ($0<q_{1}\leq2$ and non integer case) can   be  easily dealt with  by 
 the Jensen  inequality.  The theorem is then proved.  
\qed

\begin{Rem}\label{r5.3}
If the diffusion coefficient of SDE \eqref{eq5} is linear (which means $\sigma(x)=Cx$ for SDE \eqref{eq1})  then under Assumption \ref{as2r}, we have
\begin{equation*}
\EE \mid y(T)-y_{\Delta}(T)\mid^{q}\leq C\Delta^{q_{1}},
\end{equation*}
\begin{equation*}
\EE \mid x(T)-x_{\Delta}(T)\mid^{q}\leq C\Delta^{q_{1}}.
\end{equation*}
In  this case, since the diffusion coefficient is a constant, the fifth row of inequality \eqref{4aa} has a simple form:
\begin{equation*}
\EE \mid y_{\Delta}(t)-\bar{y}_{\Delta}(t)\mid^{p}\leq 2^{p-1}[\Delta^{p}\big( h(\Delta)\big)^{p}+\frac{p(p-1)}{2}\Delta^{\frac{p}{2}}],\quad \text for \quad t\in[t_{k},t_{k+1}).
\end{equation*}
Since $\Delta^{\frac{1}{4}}h(\Delta)\leq c_{0}$, it is not difficult to prove that
\begin{equation*}
\EE \mid y_{\Delta}(t)-\bar{y}_{\Delta}(t)\mid ^{p}\leq C^{p}p^{p}\Delta^{\frac{p}{2}}.
\end{equation*} 
The above conclusions follow from similar argument in the 
 proof of Theorem \ref{4.1}.
\end{Rem}

\begin{example}
Lwt return to  Example \ref{a}. Since the coefficients of \eqref{tcev} satisfy
\begin{equation*}
f'(y)=-k\lambda e^{-y}+\theta^{2}(1-\alpha)e^{-2(1-\alpha)y},~ g'(y)=-(1-\alpha)\theta e^{-(1-\alpha)y},
\end{equation*}
and $\alpha\in(\frac{3}{4},1)$, there exists a positive constant K, such that
\begin{equation*}
f'(y)+\frac{q-1}{2}g'(y)^{2}\leq K, \forall y\in\mathbb{R}.
\end{equation*}
Then for $\forall y_{1},y_{2}\in\mathrm{R}$, we have
\begin{equation*}
\begin{aligned}
&(y_{1}-y_{2})\big(f(y_{1})-f(y_{2})\big)+\frac{q-1}{2}\mid g(y_{1})-g(y_{2})\mid^{2}\\
=&\int_{0}^{1}f'(y_{1}+\zeta(y_{2}-y_{1}))\mathrm{d}\zeta\mid y_{1}-y_{2}\mid^{2}+\frac{q-1}{2}\bigg(\int_{0}^{1}g'(y_{1}+\zeta(y_{2}-y_{1}))\mathrm{d}\zeta\bigg)^{2}\mid y_{1}-y_{2}\mid^{2}\\
\leq&\int_{0}^{1}[f'(y_{1}+\zeta(y_{2}-y_{1}))+\frac{q-1}{2}\big( g'(y_{1}+\zeta(y_{2}-y_{1}))\big)^{2}]\mathrm{d}\zeta\mid y_{1}-y_{2}\mid^{2}\\
\leq& K\mid y_{1}-y_{2}\mid^{2}.
\end{aligned}
\end{equation*}
Thus   Theorem \ref{4.1} can be applied to Example  \ref{a}.
\end{example}
\begin{example}
Similarly, we can prove easily that if $r>4\rho-3$ the coefficients of \eqref{tair}  in Example \ref{b} satisfy 
\begin{equation*}
(y_{1}-y_{2})\big(f(y_{1})-f(y_{2})\big)+\frac{q-1}{2}\mid g(y_{1})-g(y_{2})\mid^{2}\leq K\mid y_{1}-y_{2}\mid^{2}, 
\end{equation*}
for all $q>2$, and for all  $y_{1},y_{2}\in\mathbb{R}$. This means that  Theorem \ref{4.1} can be applied to Example \ref{b}.
\end{example} 

\begin{Exe}\label{GLEQ}
Consider the scalar stochastic Ginzburg-Landau SDE
\begin{equation}\label{GLE}
\mathrm{d}x(t)=[-x(t)^{3}+(\lambda+\frac{\sigma^{2}}{2})x(t)]\mathrm{d}t+\sigma x(t)\mathrm{d}W(t), x(0)=x_{0}>0.
\end{equation}
The above SDE can be considered  as a special case of the Ait-Sahalia process with $(a_{-1}, a_{0}, a_{1}, a_{2}, r, \rho) = (0, 0, \lambda+\frac{\sigma^{2}}{2}, 3, 1)$, and satisfies $r>2\rho-1$. As mentioned above, this SDE has positive solution, and the boundaries are unattainable. We already know that the analytic solution of this SDE can be expressed as
\begin{equation}\label{GLEex}
x(t)=\frac{x_{0}exp(\lambda t+\sigma W(t))}{\sqrt{1+2x_{0}^{2}\int_{0}^{t}exp(2\lambda s+2\sigma W(s))\mathrm{d}s}},
\end{equation}
from which we can also see that the solution is actually positive, if $x_{0}>0$.
After using the logarithmic transformation $y(t)=\ln x(t)$, we get
\begin{equation}
\mathrm{d}y(t)=[-e^{2y(t)}+\lambda]\mathrm{d}t+\sigma\mathrm{d}W(t).
\end{equation}\label{GLEt}
Apparently, the coefficients of the above SDE satisfy Assumption \ref{as2r}, and 
then  Theorem \ref{4.1}  can be applied to this process.
\end{Exe}

The  examples mentioned above take values in $(0,+\infty)$, so the new processes obtained  by the logarithmic transformation take value in the whole line. But there are processes whose solutions take value in a certain domain $(0,D)$, and in this case we can change the logarithmic transformation to an appropriate transformation accordingly.
\begin{Exe}\label{SISEPI}
Consider the following stochastic susceptible-infected-susceptible (SIS) model (appeared in \cite{Gray2011})
\begin{equation*}
\left \{
\begin{aligned}
&\mathrm{d}S(t)=[\mu M-\beta S(t)I(t)+\gamma I(t)-\mu S(t)]\mathrm{d}t-\sigma S(t)I(t)\mathrm{d}W(t),\\
&\mathrm{d}I(t)=[\beta S(t)I(t)-(\mu+\gamma)I(t)]\mathrm{d}t+\sigma S(t)I(t)\mathrm{d}W(t),
\end{aligned}
\right.
\end{equation*}
with initial values $I(0)+S(0)=I_{0}+S_{0}=M$.
Since,
\begin{equation*}
\mathrm{d}[S(t)+I(t)]=[\mu M-\mu(S(t)+I(t)\big)]\mathrm{d}t,
\end{equation*}
 it is not difficult to get that $I(t)+S(t)=M$. Then we only need to  study the following SDE
\begin{equation}\label{SISI}
\mathrm{d}I(t)=[\big(\beta M-\mu-\gamma\big)I(t)-\beta I(t)^{2}]\mathrm{d}t+\sigma I(t)\big( M-I(t)\big)\mathrm{d}W(t), 
\end{equation}
with initial value $I(0)=I_{0}\in(0,M)$. The coefficients of SDE \eqref{SISI}  satisfy the local Lipschitz condition, then we know this SDE will have unique maximal local solution on $t\in[0,\tau_{\varepsilon})$, here the $\tau_{\varepsilon}$ is the time  to zero or $M$ for the first time.
Denote $\varsigma=\beta M-\mu-\gamma$, the scale function has the following form
\begin{equation*}
\begin{aligned}
p(x)=&\int_{c}^{x}exp\{-2\int_{c}^{\xi}\frac{\varsigma u-\beta u^{2}}{\sigma^{2}u^{2}(M-u)^{2}}\mathrm{d}u\}\mathrm{d}\xi\\
=&\int_{c}^{x}e^{\frac{-2\varsigma}{\sigma^{2}}\int_{c}^{\xi}\frac{1}{u(M-u)^{2}}\mathrm{d}u+\frac{2\beta}{\sigma^{2}}\int_{c}^{\xi}\frac{1}{(M-u)^{2}}\mathrm{du}}\mathrm{d}\xi\\
=&\int_{c}^{x}e^{\frac{-2\varsigma}{\sigma^{2}M}\big(\frac{1}{M-\xi}-\frac{1}{M-c}\big)+\frac{2\varsigma}{\sigma^{2}M^{2}}[\ln\frac{M-\xi}{\xi}-\ln\frac{M-c}{c}]+\frac{2\beta}{\sigma^{2}}\big(\frac{1}{M-\xi}-\frac{1}{M-c}\big)}\mathrm{d}\xi\\
=&C\int_{c}^{x}e^{\frac{2(\mu+\gamma)}{\sigma^{2}M}\frac{1}{M-\xi}}\big(\frac{M-\xi}{\xi}\big)^{\frac{2\varsigma}{\sigma^{2}M^{2}}}\mathrm{d}\xi,
\end{aligned}
\end{equation*}
where we used the known results 
\[
\int\frac{\mathrm{d}x}{x(a+bx)^{2}}=\frac{1}{a(a+bx)}-\frac{1}{a^{2}}\ln\frac{a+bx}{x}\,.
\]
Obviously
\begin{equation*}
p(M-)=+\infty\,. 
\end{equation*}
Now we consider $P(0+)$.  
When $\frac{2\varsigma}{\sigma^{2}M^{2}}\geq1$ it is easy to check 
\begin{equation*}
p(0+)=-\infty.
\end{equation*}
If  $\frac{2\varsigma}{\sigma^{2}M^{2}}<1$, then for any $0<x<c$
($c>0$ is a number in $(0, M)$),  we have 
\begin{equation*}
\begin{aligned}
v(x)=&\int_{x}^{c}\int_{y}^{c}p'(y)\frac{2}{p'(z)\sigma^{2}(M-z)^{2}z^{2}}\mathrm{d}z\mathrm{d}y\\
=&\int_{x}^{c}\int_{y}^{c}e^{2\int_{y}^{z}[\frac{\varsigma}{\sigma^{2}}\frac{1}{(M-u)^{2}u}-\frac{\beta}{\sigma^{2}}\frac{1}{(M-u)^{2}}]\mathrm{d}u}\frac{2}{\sigma^{2}(M-z)^{2}z^{2}}\mathrm{d}z\mathrm{d}y\\
=&\frac{2}{\sigma^{2}}\int_{x}^{c}\int_{y}^{c}e^{\frac{2(\mu+\gamma)}{\sigma^{2}M}[\frac{1}{M-y}-\frac{1}{M-z}]}
\big(\frac{M-y}{y}\big)^{\frac{2\varsigma}{\sigma^{2}M^{2}}}\big(\frac{z}{M-z}\big)^{\frac{2\varsigma}{\sigma^{2}M^{2}}-2}\frac{1}{(M-z)^{4}}\mathrm{d}z\mathrm{d}y\\
\geq& C\int_{x}^{c}\big(\frac{M-y}{y}\big)^{\frac{2\varsigma}{\sigma^{2}M^{2}}}\frac{1}{(M-y)^{2}}\int_{y}^{c}\big(\frac{z}{M-z}\big)^{\frac{2\varsigma}{\sigma^{2}M^{2}}-2}\frac{1}{(M-z)^{2}}\mathrm{d}z\mathrm{d}y\\
\geq& C\int_{x}^{c}\big(\frac{M-y}{y}\big)^{\frac{2\varsigma}{\sigma^{2}M^{2}}}\int_{y}^{c}\big(\frac{z}{M-z}\big)^{\frac{2\varsigma}{\sigma^{2}M^{2}}-2}\big(\frac{z}{(M-z)}\big)^{'}\mathrm{d}z\mathrm{d}y\\
\geq&C\big\{ \int_{x}^{c}\frac{M-y}{y}\mathrm{d}y-\int_{x}^{c}\big(\frac{M-y}{y}\big)^{\frac{2\varsigma}{\sigma^{2}M^{2}}}\mathrm{d}y\big\}.
\end{aligned}
\end{equation*}
Thus $v(0+)=+\infty$ for $\frac{2\varsigma}{\sigma^{2}M^{2}}<1$.  Hence  we  conclude  that SDE \eqref{SISI} will have an unique global solution $I(t)\in(0,M)$, and $I(t)$ will not attain the boundaries $0$ and $M$ in finite time.  This means that the disease will neither die out nor become a pandemic in finite time.  Let $ y(t)=\ln I(t)-\ln\big(M-I(t)\big)$. Then $ I(t)=M-\frac{M}{1+e^{y(t)}}$ will be always in $(0, M)$, and
  $M-I(t)=\frac{M}{1+e^{y(t)}}$ will be always in $(0, M)$ as well. 
  An application of It\^{o} formula
 with $y(t)=\ln I(t)-\ln \big(M-I(t)\big)$ yields 
 \begin{equation}\label{sist}
 \mathrm{d}y(t)=f(y(t))\mathrm{d}t+g(y(t))\mathrm{d}W(t),
 \end{equation}
 where 
 \[
 f(y)=\varsigma-(\mu+\gamma)e^{y}+\frac{\sigma^{2}M^{2}}{2}-\frac{\sigma^{2}M^{2}}{1+e^{y}}\,,\quad g(y)=\sigma M\,.
 \] 
 Since $g(y)=\sigma M$,
 we have 
\begin{equation*}
yf(y)+\frac{p-1}{2}g^{2}(y)\leq Kp, ~\forall p>2.
\end{equation*}The derivative of drift coefficients $f(y)$ is as follows
 \begin{equation*}
 f'(y)=-(\mu+\gamma)e^{y}+ \frac{\sigma^{2}M^{2}e^{y}}{(1+e^{y})^{2}}.
 \end{equation*}
 It is easy to get that
 \begin{equation*}
 \lim\limits_{y\rightarrow-\infty}f'(y)=0,~~\lim\limits_{y\rightarrow+\infty}f'(y)=-\infty.
 \end{equation*}
 Therefore there exists a positive constant such that $f'(y)\leq K$, for any $y\in\mathbb{R}$. As a consequence, we have 
 \begin{equation*}
 (y_{1}-y_{2})\big(f(y_{1})-f(y_{2})\big)\leq K(y_{1}-y_{2})^{2}.
 \end{equation*}
Moreover,   we have 
 \begin{equation*}
 \begin{aligned}
 \mid f'(y)\mid=&\mid \frac{\sigma^{2}M^{2}e^{y}}{(1+e^{y})^{2}}-(\mu+\gamma)e^{y}\mid\\
 \leq &\max\{\sigma^{2}M^{2},(\mu+\gamma)\}(e^{y}+e^{-y})\\
 \leq&Ke^{\mid y\mid}.
 \end{aligned}
 \end{equation*}
 This implies 
 \begin{equation*}
 \mid f(y_{1})-f(y_{2})\mid\leq K(e^{\mid y_{1}\mid}+e^{\mid y_{2}\mid})\mid y_{1}-y_{2}\mid.
 \end{equation*}
To summarize  we have proved  the coefficients of the SDE \eqref{sist} satisfy the Assumption \ref{as2r}. Hence we can apply all theorems in this work  to this process.  Especially, using this transformation, we can construct a numerical approximation which can maintain the domain of the original solution.
\end{Exe}

\section{Numerical exprements}
To validate  our numerical methods, several numerical experiments are given in this section. 

First, the GLE model \eqref{GLE} in Example \ref{GLEQ} is used as an example to illustrate the positivity and the convergence of our numerical method. The GLE model has an explicit  expression \eqref{GLEex} for the exact solution, but there is an integral term in the expression, which can not be accurately evaluated. We  approximate the integral by a Riemann sum with step size $\Delta=2^{-12}$. For this model due to the availability of this explicit expression for the exact solution, we can compare the exact solution and the approximate solution by our approach.

The  truncated Euler-Maruyama method (TE) proposed in \cite{Mao2015} was proved to be strongly convergent under the coercive condition. The coefficients of GLE model \eqref{GLE} satisfy the coercive condition, so numerical approximation obtained by using TE method will converge to the exact solution $x(t)$ of SDE \eqref{GLE}.  However, the positivity  of this method cannot be guaranteed. The idea of logTE method proposed in this paper is to apply the TE method to SDE \eqref{GLEt} rather than SDE \eqref{GLE} directly. Then transforming back yields a positive approximation $x_{\Delta} (t)=e^{y_{\Delta}(t)}$ for the original solution $x(t)$ of SDE \eqref{GLE}.

\begin{figure}[h]
\centering
\subfigure[Sample path]{\includegraphics[width=0.45\textwidth]{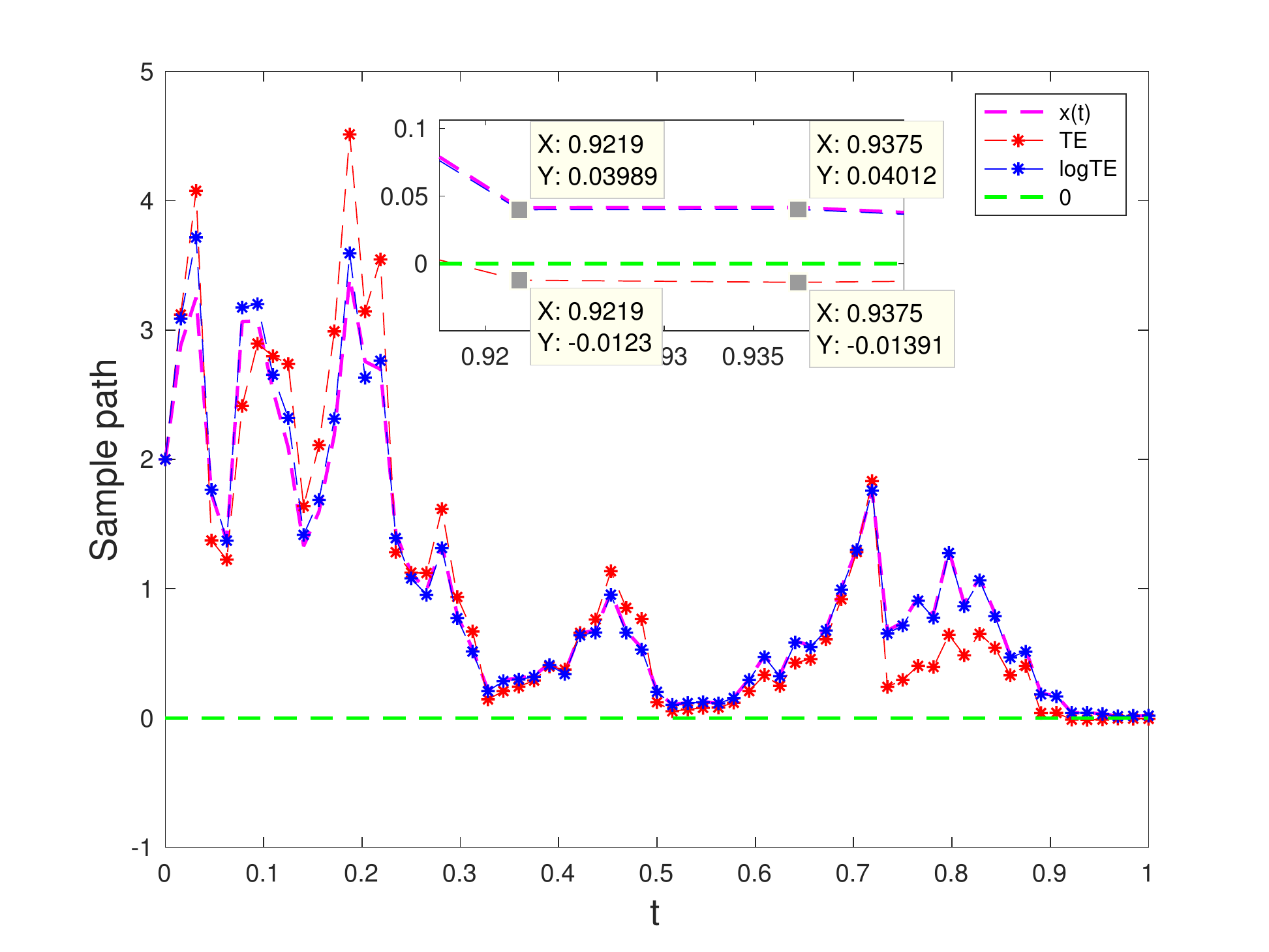}}
\subfigure[Convergence rate]{\includegraphics[width=0.45\textwidth]{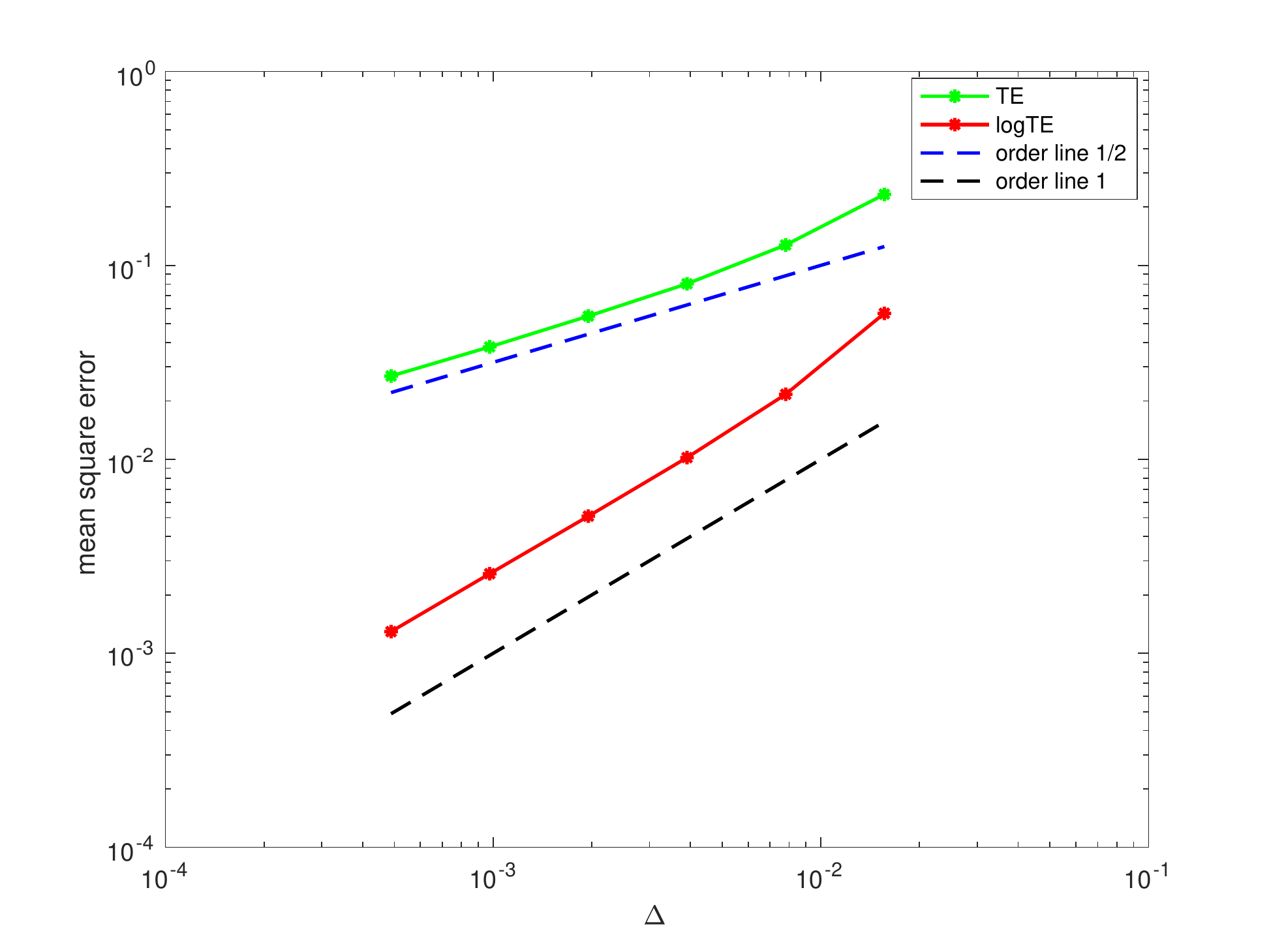}}
 \caption{Simulations of  GLE \eqref{GLEQ}.}
 \label{fig1}
\end{figure}

In our experiments,  \eqref{GLE} is numerically solved by TE method and logTE method respectively. The parameters we choose are $\lambda=1,\sigma=4,T=1, x_{0}=2$.  
We plot the single trajectories of the exact solution and two numerical solutions with equidistant step size $h = 2^{-6}$ in (a) of Figure \ref{fig1}. The results show that compared with TE method, our method does have advantages in preservation of positivity.
 
 We compute six different numerical solutions  with step sizes $\Delta= 2^{-6}, 2^{-7},\cdots, 2^{-11}$ of  the logTE method and the TE method on the same Brownian path,  and gain the corresponding mean square errors. The mean square  error at the terminal time $T=1$, i.e., $\| x(T)-x_{\Delta}(T)\|_{L_{2}}$, is given by 
 \begin{equation*}
\| x(T)-x_{\Delta}(T)\|_{L_{2}}= \bigg(\frac{1} { M_{0}}\sum_{i=1}^{M_{0}} \mid x^{i}(T)-x^{i}_{\Delta}(T) \mid^{2}\bigg)^{\frac{1}{2}}, 
 \end{equation*}
 where $M_{0}=5000$,  $x^{i}(T)$ and $x^{i}_{\Delta}(T)$ denote the number of sample paths, the ith exact solution and the ith numerical solution respectively. We show mean square errors of logTE method and TE method in (b) of Figure \ref{fig1}. Since the diffusion coefficient of the transformed SDE \eqref{GLEt} is a constant, then according to Remark \ref{r5.3} we know that the convergence rate of logTE method is 1 in this case. Since the convergence order of TE method is close to 0.5, that is to say, for the GLE model, our method is superior to the TE method in terms of convergence rate. The result showed in (b) of Figure \ref{fig1} is consistent with our theoretical results

\begin{figure}[h]
\centering
\subfigure[Sample path]{\includegraphics[width=0.45\textwidth]{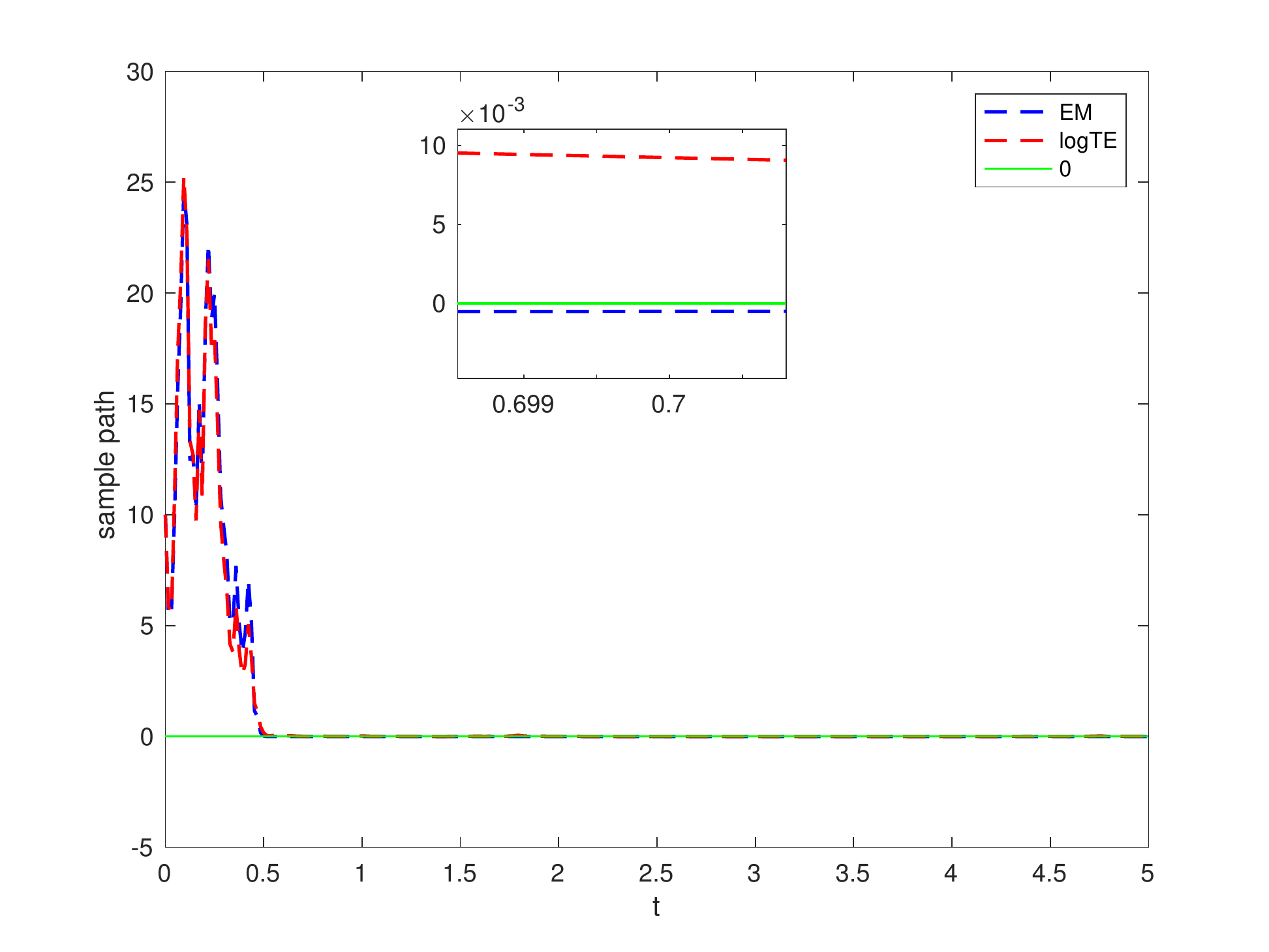}}
\subfigure[Convergence rate]{\includegraphics[width=0.45\textwidth]{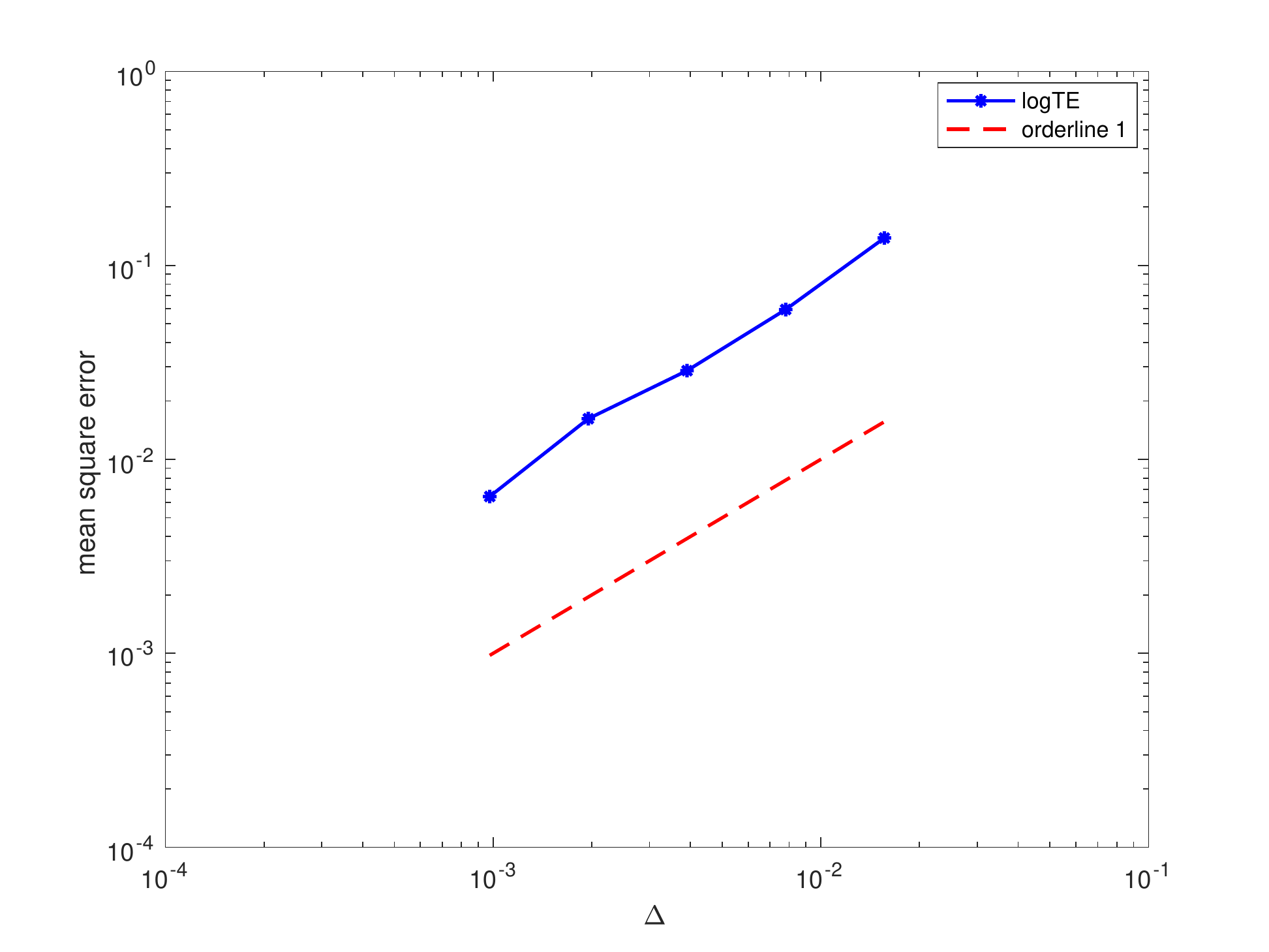}}
 \caption{Simulations of SIS model  \eqref{SISI}.}
 \label{fig2}
\end{figure}

Next, we consider the SIS model \eqref{SISI}. The parameters in this model are set  as follows:  $\beta=0.5, N=100, \mu=20, \gamma=25,\sigma=0.035,I_{0}=10 $.  Again, the sample path and  convergence rate of the logTE method are shown in the Figure 2. For comparison, we plot single trajectories of the  Euler-Maruyama method  and  logTE method with equidistance $\Delta=2^{-6}$  in (a) of Figure \ref{fig2}, from which we can see that approximation from our method can maintain the domain of the exact solution while the Euler-Maruyama approximation cannot. Since the explicit solution of this SDE can not be derived, the numerical solution with step size $\Delta=2^{-14}$ is used as the reference solution. We compute five different numerical solutions of the logTE with step size $\Delta=2^{-6},, 2^{-7}, \cdots,2^{-10}$, and show the estimated mean square convergence rate in (b) of Figure \ref{fig2}. 
 From Figure \ref{fig2}, we can also find that the numerical results are consistent with our theoretical results given above.

\begin{figure}[h]
\centering
\subfigure[Sample path]{\includegraphics[width=0.45\textwidth]{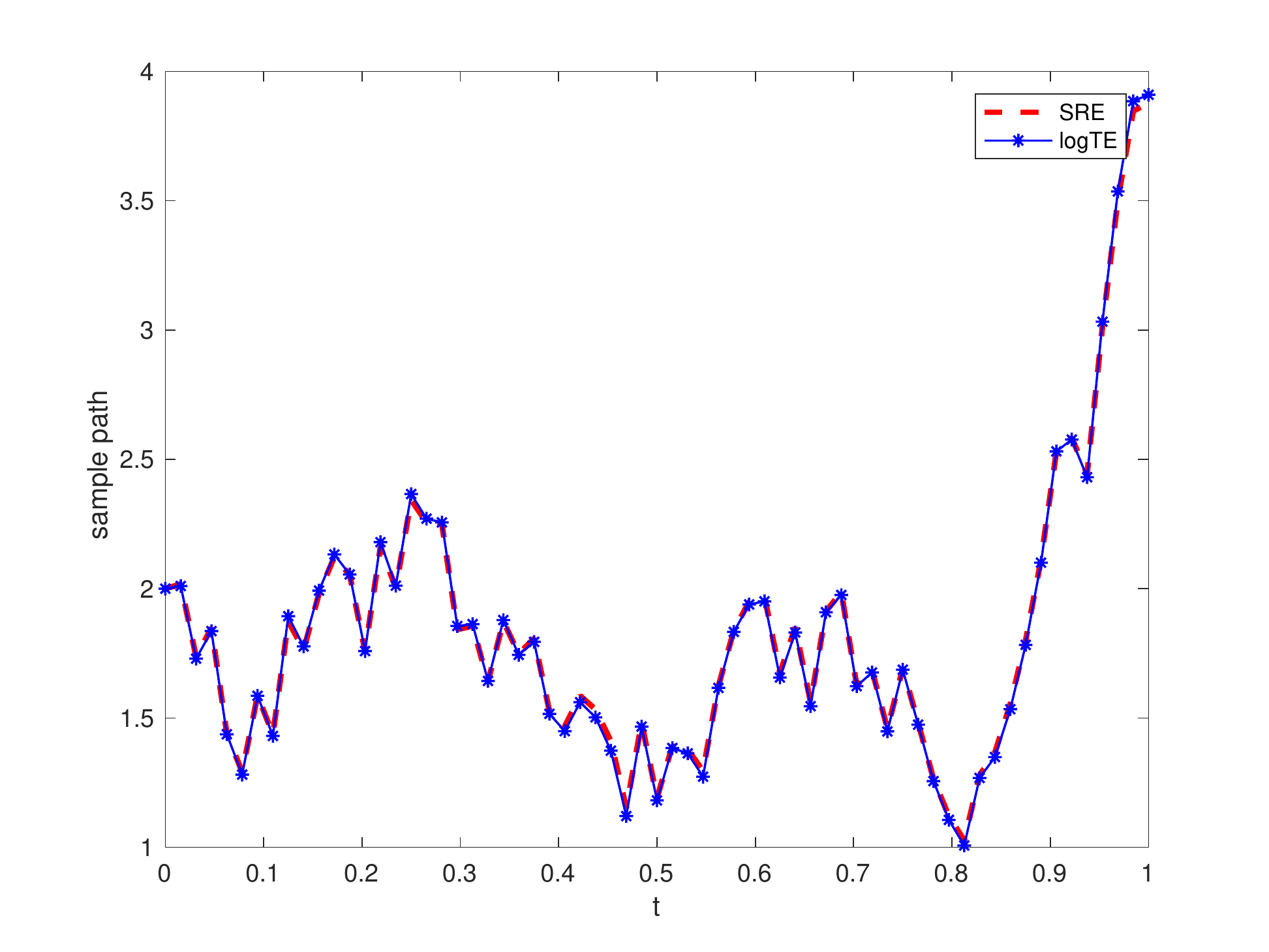}}
\subfigure[Convergence rate]{\includegraphics[width=0.45\textwidth]{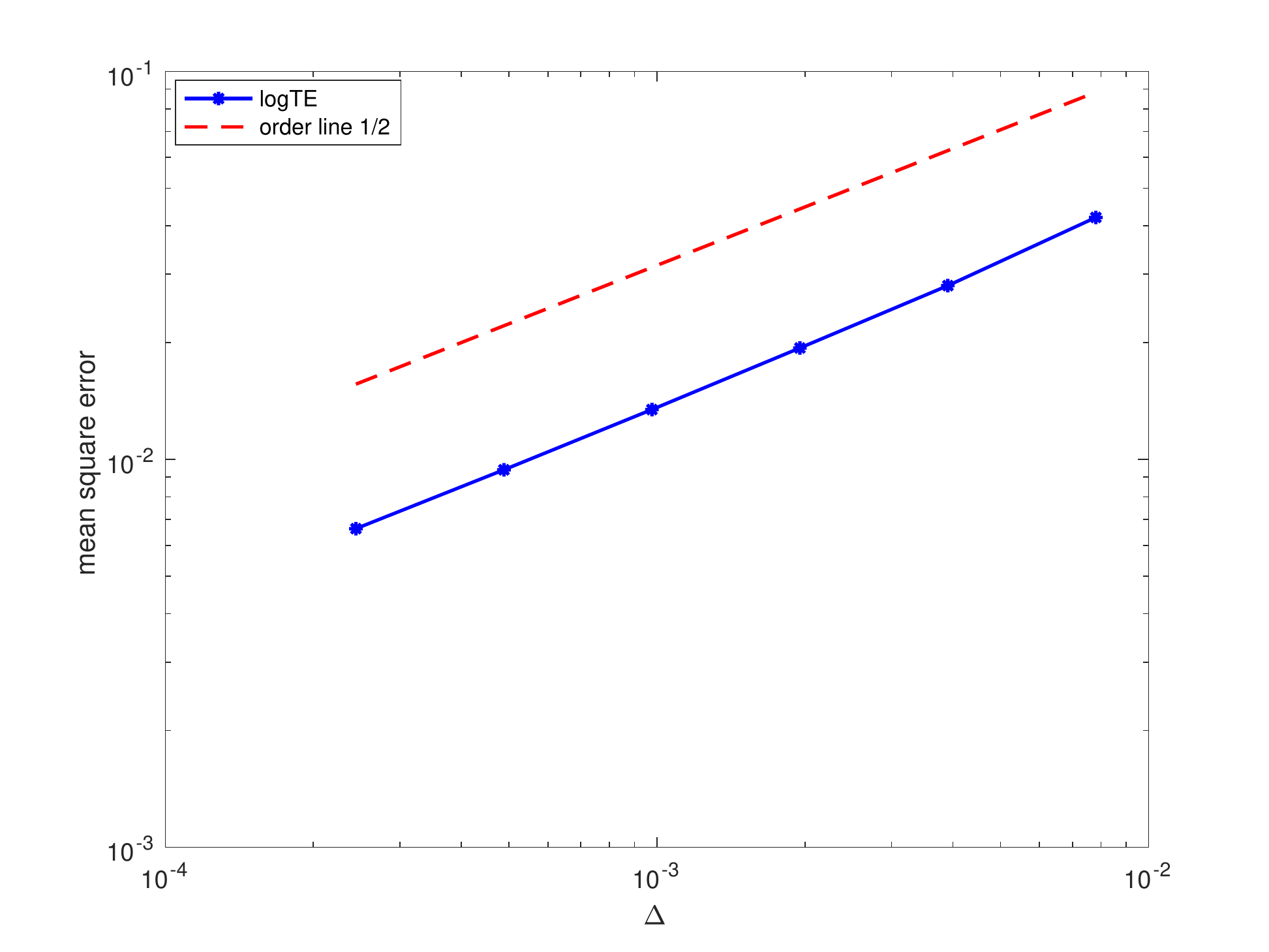}}
 \caption{Simulations of CIR model \eqref{ciro}.}
 \label{fig3}
\end{figure}

Besides, we also consider the CIR model \eqref{ciro}. Set $\kappa=2, \lambda=1,\theta=\frac{3}{2},x_{0}=2$. Since the exact solution has no explicit expression, we use the drift-implicit square-root Euler method (SRE) proposed in \cite{Dereich2012} which can preserve the positiveness as a reference solution. We draw single trajectories of SRE with $\Delta=2^{-14}$ and logTE with $\Delta=2^{-6}$ in (a) of Figure \ref{fig3}. In addition, the coefficient of the SDE \eqref{cirt} does not satisfy Assumption \ref{as2r}, so we do not get the convergence rate of the logTE method for this model. We show the convergence order of this method in (b) of Figure \ref{fig3} by six different numerical solutions  with step size $2^{-7}, 2^{-8}, \cdots,2^{-12}$.

\begin{figure}[h]
\centering
\subfigure[Sample path]{\includegraphics[width=0.45\textwidth]{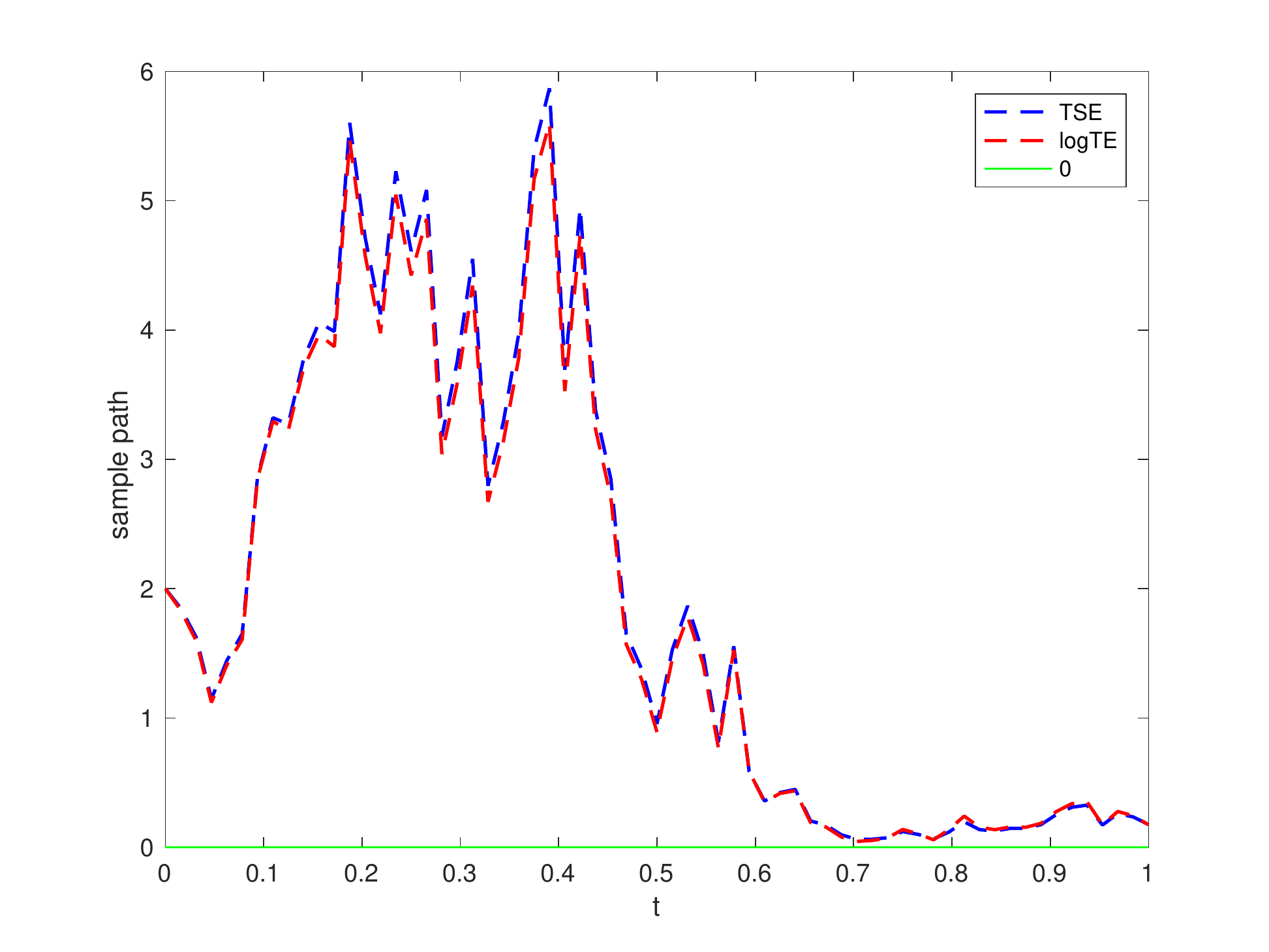}}
\subfigure[Convergence rate]{\includegraphics[width=0.45\textwidth]{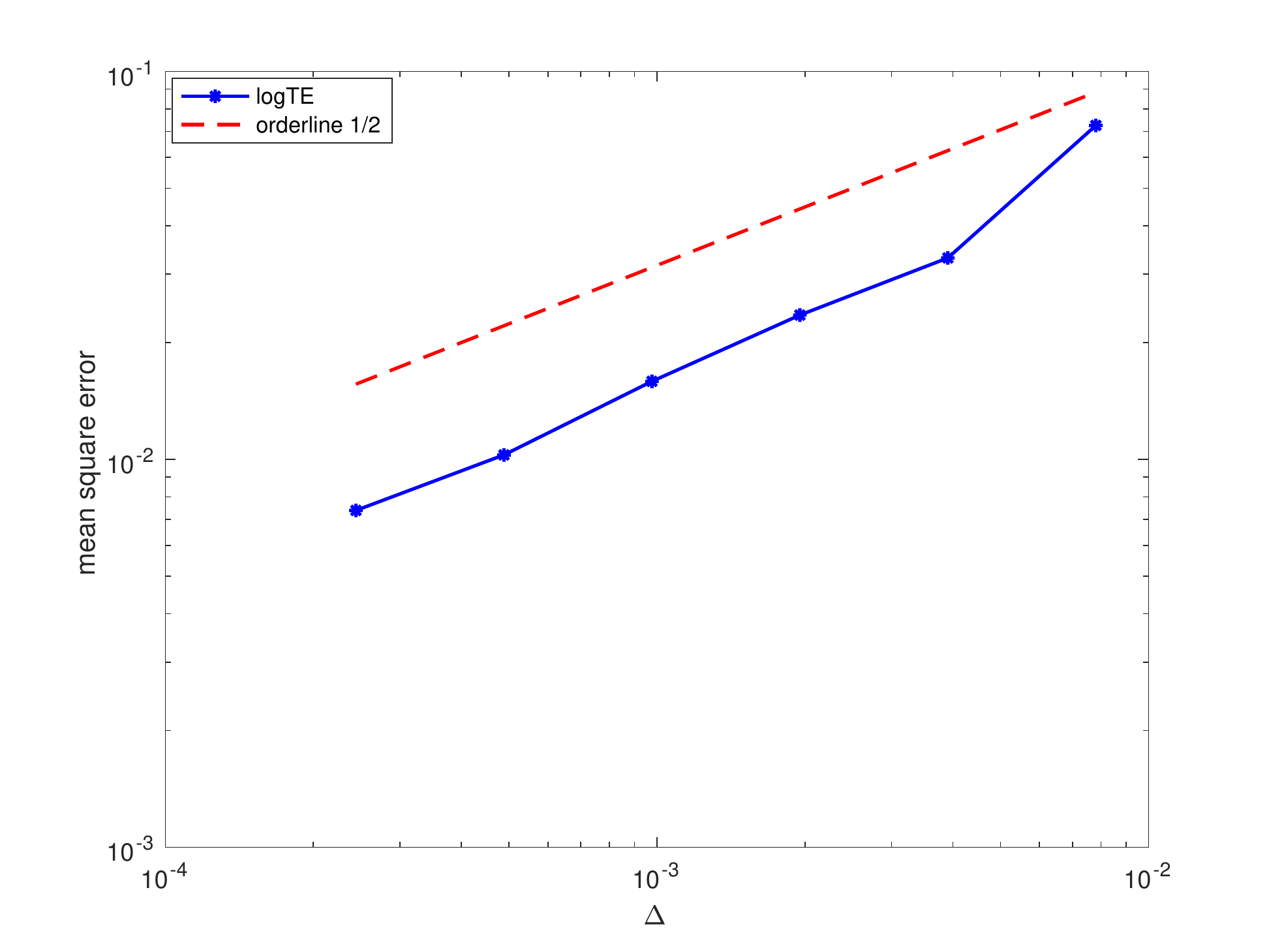}}
 \caption{Simulations of CEV model  \eqref{cev}.}
 \label{fig4}
\end{figure}

For the last example we consider the CEV process \eqref{cev}. Set $\kappa=4, \lambda=0.5, \theta=3, x_{0}=2, \alpha=\frac{7}{8}$. Similarly, since the exact solution of this model is difficult to get, we will use the approximation from semi-implicit method (denote as TSE)  appeared in \cite{Neuenkirch2014} as a reference solution. We compare  single trajectories of  the logTE method  and the TSE method with same step size $\Delta=2^{-6}$ in (a) of Figure \ref{fig4}.  
The numerical solution from TSE was proved to be able to maintain the positiveness of the original solution. From (a) of Figure \ref{fig4}, we can see that the trajectory of approximation obtained by logTE method is very close to that of TSE.  The estimated mean square errors of logTE numerical solutions with six different step sizes $\Delta=2^{-7},2^{-8}, \cdots, 2^{-12}$ are shown in (b) of Figure \ref{fig4}, from which one clearly observes that mean square convergence rate is $\frac{1}{2}$ for logTE.

\section{Conclusion}

In this paper, we propose a class of logarithmic Euler-Maruyama type methods for SDEs with positive solutions. Our idea is to first apply the logarithmic transformation $y(t)=\ln x(t)$, and then obtain a numerical approximation $y_{\Delta}(t)$ to the transformed process $y(t)$, $t\in[0,T]$ by using Euler-Maruyama type method. And $x_{\Delta}(t)=e^{y_{\Delta}(t)}$ will be a positive approximation to the original solution $x(t)$.  We demonstrate some convergence and convergence rate conclusions for these methods in path-wise  and $L^{p}$ sense.  The using of logarithmic transformation makes the $L^{p}$ convergence analysis much more difficult,
since the coefficients of the transformed process $y(t)$ may grow exponentially and we need to deal with the exponential integrability of both $y(t)$ and its approximated solution imposed by our transformation. 
We also note that for SDE whose solution takes value in a certain domain $(0,D)$, approximate solution from our method is inside the domain of the original solution if the logarithmic transformation is changed accordingly, see Example \ref{SISEPI}.

\section*{Acknowledgements}
Y.Yi and J. Zhao are supported by the National Natural Science Foundation of China (11771112, 11671112).
Y. Hu is supported by   an NSERC discovery grant   and a startup fund from   University of Alberta at Edmonton.

\section*{References}

\end{document}